\documentclass[11pt]{article}
\usepackage{amsfonts, amsmath, amscd,amssymb}%,fullpage}
\usepackage{amssymb, amsbsy, amsthm}%,  amstext}%, amsopn, verbatim,multicol}
\usepackage{latexcad}
\newfont{\eufm}{eufm10 scaled\magstep1}

\newtheorem{thm}{Theorem}[section]
\newtheorem{lem}[thm]{Lemma}
\newtheorem{cor}[thm]{Corollary}

\numberwithin{equation}{section}

\newtheorem*{klem}{Key Lemma}

\newcommand{\LL}{\mathcal{L}}

\newcommand{\OO}{\mathcal{O}}

\newcommand{\NN}{\mathbb{N}}
\newcommand{\Z}{\mathbb{Z}}

\newcommand{\C}{\mathbb{C}}

 \newcommand{\ttt}{\textsf}

\newcommand{\lgr}{\ttt{\text{-}grmod}}
 \newcommand{\lGr}{\ttt{\text{-}Grmod}}

 \newcommand{\ltors}{\text{-}\ttt{tors}}
 \newcommand{\lTors}{\text{-}\ttt{Tors}}

 \DeclareMathOperator{\spc}{Spec}

\DeclareMathOperator{\coh}{\ttt{Coh} }
\DeclareMathOperator{\Qcoh}{\ttt{Qcoh} }

\begin{document}
\title{
%An Equivalence of Categories for
Noncommutative Deformations of Type A Kleinian Singularities and
Hilbert Schemes.}

\author{Ian M. Musson\thanks{partially supported by NSF grant  DMS-0099923.}
\\ Department of Mathematical Sciences\\
University of Wisconsin-Milwaukee, USA\\E-mail: musson@uwm.edu}
 \maketitle

\begin{abstract}
Let $H_{\mathbf{k}}$ be a symplectic reflection algebra
corresponding to a cyclic subgroup $\Gamma \subseteq SL_2 \C$ of
order $n$ and   $U_{\mathbf{k}} = eH_{\mathbf{k}} e$ the spherical
subalgebra of
   $H_{\mathbf{k}}$.
 We show that for suitable ${\mathbf{k}}$ there is a filtered
 $\Z^{n-1}$-algebra $R$ such that
 \begin{itemize}
   \item[{(1)}] there is an equivalence of categories
   $R-\mathrm{qgr} \simeq U_{\bf k}$-mod ,
   \item[{(2)}]  there is an equivalence of categories $gr R-\mathrm{qgr} \simeq
   \coh(Hilb_\Gamma
\mathbb{C}^2)$.

 \end{itemize}
Here $ \coh(Hilb_\Gamma \mathbb{C}^2)$ is the category of coherent
sheaves on the $\Gamma$-Hilbert scheme. and for a graded algebra
$\mathcal{R},$ we write $ \mathcal{R}-\mathrm{qgr}$ for the
quotient category of
finitely generated graded $\mathcal{R}$-modules modulo torsion. \\
\end{abstract}

\section{Introduction}\label{intro}
 This paper is motivated by a construction of Iain Gordon
and Toby Stafford concerning the representation theory of a
symplectic reflection algebra $H_c$ and the
   spherical subalgebra $U_c$ of $H_c.$ This construction addresses  issues
   raised in Conjecture 1.6 in [GK].
 We refer to [EG] for the
   definition of symplectic
reflection algebra  and to [GS1] for background and the analogy
for   with the Beilinson-Bernstein theorem.\\
 \indent Suppose that $G$ is a finite subgroup of $SL_2
\mathbb{C}$, and let $\Gamma = G \wr S_m$, the wreath product of
$G$ by the symmetric group $S_m$. The group $\Gamma$ acts on $V =
(\mathbb{C}^2)^m$ preserving the natural sympletic structure. Let
$Y_{\Gamma,m}$ be the set of $\Gamma$-invariant ideals $I$ in the
Hilbert scheme of $m|G|$ points in
 $\mathbb{C}^2$ such that the quotient $\mathbb{C}[x,y]/I$ is isomorphic to a direct sum of
 $m$ copies of the regular representation of $G.$
By [W, Theorem 4.2] there is a crepant resolution
\[ Y_{\Gamma,m}  \longrightarrow V/\Gamma .\]
This gives the bottom arrow in the diagram below. The algebra
$U_c$ has a filtration such that the associated graded algebra is
isomorphic to $\mathcal{O}({V/\Gamma})\mathrm{-mod}.$ The passage
to associated graded modules is indicated by the vertical arrow on
the left.  Gordon and Stafford suggest that there is a suitable
category
 that will complete the diagram.
\begin{center}
\begin{picture}(120,76)
\thinlines \drawvector
 {51.0}{62.0}{62.0}{1}{0}
\drawvector{51.0}{14.0}{32.0}{1}{0}
\drawcenteredtext{18.0}{62.0}{$U_c{-mod}$}
\drawlefttext{125.0}{62.0}{$???$}
\drawcenteredtext{18.0}{14.0}{$\mathcal{O}({V/\Gamma})\mathrm{-mod}$}

\drawlefttext{95.0}{14.0}{$\mathcal{O}_{Y_{\Gamma,m} }
\mathrm{-mod } $}
   \drawvector{132.0}{54.0}{32.0}{0}{-1}

\drawvector{18.0}{54.0}{32.0}{0}{-1}
\drawcenteredtext{67.0}{70.0}{$\simeq$}
%\drawcenteredtext{67.0}{6.0}{$\simeq$}
\drawcenteredtext{8.0}{38.0}{$gr$}
%\drawcenteredtext{78.0}{38.0}{$gr$}
\drawcenteredtext{122.0}{38.0}{$gr$}
\end{picture}
\end{center}
  \indent
The main theorem in [GS1] shows that this is possible in the
crucial special case where $G = 1$, and so $\Gamma = S_m.$
Applications of this result are given in [GS2].  Another special
case arises
   when $m= 1,$ and so $G = \Gamma.$  Then $V/\Gamma$ is a Kleinian singularity and the algebras $H_c$
   were actually introduced earlier by Hodges [H] if $\Gamma$ is cyclic, and by
    Crawley-Boevey and Holland [CBH] in general.
    The purpose of
   this paper is to solve the problem of Gordon and Stafford when $\Gamma$ is cyclic of
   order $n.$ Note that when $\Gamma$ is the symmetric group or a finite subgroup of $SL_2
\mathbb{C}$ then $Y_{\Gamma,m}$ is the same as the
$\Gamma$-Hilbert scheme $Hilb_\Gamma\mathbb{C}^2.$

   To state the main theorem requires some notation.  There is an action of
$\Gamma$  on the first Weyl algebra $\mathbb{C}[\partial, y]$ and
on the localization $\mathbb{C}[\partial, y^{\pm 1}].$  For ${\bf
k} \in \mathbb{C}^{n-1},$ we construct $U_{\bf k}$ and $H_{\bf k}$
as subalgebras of the crossed product $\mathbb{C}[\partial, y^{\pm
1}]*\Gamma.$ {Our parameterization of these algebras is different
from what is usually used, hence the change in notation.} Then for
suitable ${\bf k'}, {\bf k} \in \mathbb{C}^{n-1},$ we find $U_{\bf
k'}-U_{{\bf k}}$-bimodules $B({\bf k'},{\bf k})$ and give a
sufficient condition for these bimodules to induce a Morita
equivalence. Next following [GS1], we assemble the bimodules
$B({\bf k'},{\bf k})$ to form a Morita $\mathbb{Z}^{n-1}$-algebra
$R.$ This is an algebra, without identity which is  graded by
$\mathbb{Z}^{n-1} \times \mathbb{Z}^{n-1}.$  Since the algebras
$U_{\bf k}$ and the bimodules $B({\bf k'},{\bf k})$ are contained
in $\mathbb{C}[\partial, y^{\pm 1}]*\Gamma,$ they have a
filtration given by the order of the differential operators. The
algebra $R$ inherits this filtration and we write  $gr R$ for the
associated graded algebra. The associated graded algebra of
$U_{\bf k}$ is isomorphic to $\mathcal{O}(V/\Gamma)$. We write $
\coh(Hilb_\Gamma \mathbb{C}^2)$ for the category of coherent
sheaves on $Hilb_\Gamma \mathbb{C}^2.$ For a graded algebra
$\mathcal{R},$ we write $ \mathcal{R}-\mathrm{qgr}$ for the
quotient category of finitely generated graded
$\mathcal{R}$-modules modulo torsion. This notation is explained
more fully at the start of Section 3. We can now state our
main result. The definition of dominance is immediately given before Theorem \ref{ad}.\\
\\
\noindent {\bf Main Theorem.}  If ${\bf k}$ is  dominant then\\
 (1)  there is an equivalence of categories
   $R-\mathrm{qgr} \simeq U_{\bf k}$-mod ,\\
(2) there is an equivalence of categories $gr R-\mathrm{qgr}
\simeq
   \coh(Hilb_\Gamma
\mathbb{C}^2)$.\\

A brief outline of the proof is as follows.  In the next section
we  recall the construction of  $\coh(Hilb_\Gamma \mathbb{C}^2)$
as a toric variety from [IN]. We also show that there is a
$\mathbb{N}^{n-1}$ -graded ring $S$ such that the category $S
\mathrm{-qgr}$ is equivalent to $\coh(Hilb_\Gamma \mathbb{C}^2).$
We construct a $\mathbb{Z}^{n-1}$-algebra $\widehat{S}$ such that
the categories $ \widehat{S} \mathrm{-qgr}$ and $S \mathrm{-qgr}$
are equivalent.  The algebra $gr R$ is also
$\mathbb{N}^{n-1}$-graded and it is not hard to show that
$\widehat{S} \subseteq gr R$.  The proof of (2) is completed using
a Poincar\'{e} series argument as in [GS1, Section 6] to show that
$\widehat{S} = gr R$. The proof of (1) uses a generalization of
the $\mathbb{Z}$-algebra machinery developed in [GS1, Section 5].

 I would like to thank Iain Gordon for
suggesting the problem considered in this paper, and for outlining
the approach which led to a solution. The material in Section 3
and the proof of Lemma \ref{Mor} are due to him.  %We understand that
I also thank Mitya Boyarchenko of the University of Chicago for
some helpful correspondence.  Mitya has informed me that
   % made some initial progress on
   he is close to a solution of this problem for general
   Kleinian singularities.\\

\section{Toric Varieties}\label{tv}
 \indent Let $\Gamma$ be a finite subgroup of $SL_2\mathbb{C}$.
 The $\Gamma$-Hilbert scheme $Hilb_\Gamma
 \mathbb{C}^2$ is the scheme  which parameterizes $\Gamma$-invariant ideals $I$
 of $C[u,v]$ such that $\mathbb{C}[u,v]/I$ is isomorphic to the
 regular representation of $\Gamma$.  The Hilbert-Chow morphism
  \[ Hilb_\Gamma \mathbb{C}^2 \longrightarrow \mathbb{C}^2 /
  \Gamma \]
  gives a minimal  (equivalently crepant) resolution of the quotient singularity,
  [W].  For $\Gamma$ cyclic of order $n$ we construct this morphism using toric
  varieties.

 Let $N = \mathbb{Z}^2, M = Hom(N,\mathbb{Z}),$ and write $
< \;\; , \;
> : M \times N\longrightarrow \mathbb{Z}.$
 for the natural bilinear pairing. Set $v_i = (1,i)$
for $0 \leq i \leq n.$ Then let $\Delta$ be the fan in
$N_\mathbb{R} = N \otimes_\mathbb{Z} \mathbb{R}$ with one
dimensional cones $\mathbb{R}^+v_i$ and 2 dimensional cones
$\sigma_i = \mathbb{R}^+v_{i-1} + \mathbb{R}^+v_i$ for $1 \leq i
\leq n$. The cone $\sigma_i^\vee$ in $M_\mathbb{R}$ dual to
$\sigma_i$ is $\sigma_i^\vee = \mathbb{R}^+(i, -1) + \mathbb{R}^+
(1 - i, 1)$ and we set
\[A_i = \mathbb{C}[M \cup \sigma_i^\vee], \quad X_i =
Spec A_i.\] Let $X = X(\Delta) = \cup _{i = 1}^{n}
 X_i$ be the toric variety
determined by $\Delta$ and let $T = T_N$ be the dense torus acting
on $X.$  Since $det \left( \begin{array}{cc} i & -1\\1-i & 1
\end{array}\right) = 1,$ it follows from [F, Proposition, page 29] that $X$
is nonsingular.  We write elements of $A_i$ multiplicatively by
setting
\[x = (1, 0), \quad z = (0, 1) .\] Then $A_i =
\mathbb{C}[x^iz^{-1}, x^{1-i}z ]$ and the maximal ideal ${\bf
m}_i$ of $A_i$ corresponding to the $T$ fixed point $p_i \in X_i$
is
\begin{equation} \label{m}
{\bf m}_i = (x^iz^{-1}, x^{1-i}z).
\end{equation}
In order to relate $X$ to the singularity it is convenient to
introduce new indeterminates $u, v$ satisfying
\[v^n = z, \quad uv = x,\]
so that $u^n = x^{n}z^{-1}.$  Then
 \[ A_i = \mathbb{C}[u^iv^{i-n},v^{n+ 1-i}u^{1-i}] .\]
  Note that the fan
$\Delta$ is obtained by subdividing the fan $\Delta'$ with the
single  2-dimensional cone $\sigma= \mathbb{R}^+(1, 0) +
\mathbb{R}^+ (1, n).$  Also $\sigma^\vee = \mathbb{R}^+(0, 1) +
\mathbb{R}^+ (n, -1),$ so the toric variety $X'$ corresponding to
$\Delta'$ is $Spec \; C^\Gamma$ where $C =\mathbb{C}[u, v].$ It
follows easily that the map $X \longrightarrow X'$ is the minimal
resolution of the singularity.

 \indent For $1 \leq i \leq n$, let $U_i$ be
the set of $\Gamma$-invariant ideals of $\mathbb{C}[u,v]$ such
that the elements
\[ 1, u, \ldots, u^{i-1}, v, v^2, \ldots, v^{n-i} \]
form a basis for the factor algebra $\mathbb{C}[u,v]/I$.  If $I$
is a $\Gamma$-invariant ideal of $\mathbb{C}[u,v]$ such that the
factor algebra is isomorphic to $\mathbb{C}\Gamma$, then $u^n, uv,
v^n \in I$, and it follows that $I \in U_i$ for some $i$.
Furthermore if $I \in U_i$ then for a unique $(a,b) \in
\mathbb{C}^2$ we have
\[ u^i \equiv av^{n-i}, \; v^{n-i+1} \equiv bu^{i-1} \mod I. \]
This identifies $U_i$ with $\mathbb{C}^2$ and $\mathcal{O}(U_i)$
with $\mathbb{C}[u^i v^{i-n}, v^{n+1-i} u^{1-i}] = A_i$. The open
sets $U_i$ are glued together in the same way as the $X_i$ and it
follows easily that the resolution $X \longrightarrow X'$
constructed above is the Hilbert-Chow morphism, [IN].  This map is
equivariant for the action of a dense torus $T = \mathbb{C}
^\times \times \mathbb{C}^\times$. Let $\mathbb{X}(T) =
\mathbb{Z}\chi_1 \oplus \mathbb{Z} \chi_2$ be the character group
of $T.$ If $T$ acts rationally on a vector space $V$ we define the
weight space decomposition of $V$ to be $V = \bigoplus_{\chi \in
\mathbb{X}(T)} V({\chi}),$ where $V({\chi}) = \{ v \in V| \tau
\cdot v = \chi(\tau)v \; \mbox{for all} \; \tau \in T\}. $ If $dim
V({\chi})< \infty$ for all $\chi \in \mathbb{X}(T)$ we define the
Poincar\'{e} series $\mathcal{H}_V (q,t)$ of $V$
 by
\[ \mathcal{H}_V(q,t) = \sum_{r,s} \dim V({r\chi_1 + s\chi_2})\;q^r t^s.\]
We assume that  $T$ acts on $C$ so that
%for $\tau \in T$ we have$\tau \cdot u = q(\tau)u, \tau \cdot v = t(\tau)v$. Sometimes we
%abuse notation and write simply $\tau \cdot u = qu, \tau \cdot v =tv$.
$u \in C({\chi}_1)$  and $v \in C({\chi}_2).$  We have

\begin{equation} \label{la6}
x  \in C({\chi}_1 + {\chi}_2) ; \quad z  \in C(n{\chi}_2).
\end{equation}
 Let $D_{i}$ be the divisor on $X$
corresponding to $v_{i}.$ By [F, 3.4] there is
 an exact sequence.
\[ 0 \longrightarrow M \stackrel{\alpha}{\longrightarrow} \oplus_{i = 0}^{n}
\mathbb{Z}D_i \stackrel{\beta}{\longrightarrow} Pic(X)
\longrightarrow 0 . \]
 where $\alpha(m) = \sum_{i = 0}^{n} <m ,
 v_i > D_i$ and  $\beta $ sends a divisor $D$ to its class $[D]$ in $Pic(X).$ It
follows that $Pic(X)$ is generated by the $D_i$ which are subject
to the relations
\begin{equation}
 \sum_{i = 0}^{n} D_i  \; = \;  \sum_{i = 0}^{n} i D_i \; = \; 0.
\end{equation}
This implies that $ Pic(X) = \oplus_{i = 1}^{n-1} \mathbb{Z} D(i)
,$ where for $1 \leq i \leq n-1,$ we define

\begin{equation}
D(i) = \sum_{j = 0}^{i-1}   (i-j) D_{n-j}.
\end{equation}
For ${\bf b} \in \mathbb{Z}^{n-1}$ set
\begin{equation} \label{eq5}
D({\bf b}) = \sum_{i = 1}^{n-1} b_i D(i).
\end{equation}
The polynomial ring $ \mathbb{C}[x_0, \ldots, x_n]$ is graded by
$Pic(X)$ where
 we define
 \[ \deg (x_i) = [D_i] \in Pic(X) .\]
We recall  Proposition 1.1 in [Cox]. Let
\[S_{\bf b} =
 span \{s \in \mathbb{C}[x_0, \ldots, x_n]| deg(s) = [D({\bf b})] \} .\]
Then there is an isomorphism $ S_{\bf b}  \cong H^0(X,
\mathcal{O}(D({\bf b})))$
\[\phi_{\bf b}: S_{\bf b} \longrightarrow  H^0(X,
\mathcal{O}(D({\bf b})))\]  and a commutative diagram

\begin{center}
\begin{picture}(100,76)
\thinlines
%horiz vects
%top
\drawvector{-10.0}{62.0}{160.0}{1}{0}
\drawvector{56.0}{14.0}{44.0}{1}{0}
\drawcenteredtext{-44.0}{62.0}{$S_{\bf b} \otimes S_{\bf c}$}
\drawcenteredtext{164.0}{62.0}{$S_{{\bf b}+ {\bf c}}$}
\drawcenteredtext{-44.0}{14.0}{$H^0(X, \mathcal{O}(D({\bf b})))
\otimes H^0(X, \mathcal{O}(D({\bf c})))$}
\drawcenteredtext{164.0}{14.0}{$H^0(X, \mathcal{O}(D({\bf b} +
{\bf c})))$} \drawvector{164.0}{54.0}{32.0}{0}{-1}
\drawvector{-44.0}{54.0}{32.0}{0}{-1}
\end{picture}
\end{center}
%\vspace{0.5cm}
 where the horizontal arrows are multiplication and
the vertical arrows are $\phi_{\bf b} \otimes \phi_{\bf c}$ and
$\phi_{{\bf b} + {\bf c}}.$

\begin{lem} \label{NL4}
If ${\bf b} \in \mathbb{N}^{n-1},$ and $D = D({\bf b}),$ then\\
(1)  $\mathcal{O}(D({\bf b}))$ is generated by its global sections\\
(2) $H^i(X,\mathcal{O}(D)) = 0$ for $i >0.$\\
\end{lem}
\noindent {\bf Proof.} (1) Let $\psi_D$ be the piecewise linear
function associated to $D = D({\bf b})$ in  [F, Section 3.4]. It
is easy to see that if ${\bf b} \in \mathbb{N}^{n-1},$ then
$\psi_D$ is convex. The
result follows from [F, Proposition, page 68]. \\
(2) This follows from (1) and  [F, Corollary, page 74].\\

\begin{cor}  \label{NL3}
For ${\bf b}, {\bf c} \in \mathbb{N}^{n-1}$ we have
\[H^0(X,\mathcal{O}(D({\bf b})))H^0(X,\mathcal{O}(D({\bf c})))
=H^0(X,\mathcal{O}(D({\bf b} + {\bf c}))).\]
\end{cor}
\noindent
 {\bf Proof.} This follows by part (1) of the Lemma and
the first exercise on page 69 of [F].\\

Our goal is to compute a graded Poincar\'{e} series for
$H^0(X,\mathcal{O}(D)).$ To do this we need the
Atiyah-Bott-Lefschetz theorem, which we state next.
\begin{thm} \label{ABLT}
 Let $T =
\mathbb{C}^\times \times \mathbb{C}^\times$ be a 2 dimensional
torus acting on a smooth surface $X$ such that the fixed point set
$X^T$ is finite.  Let $\mathcal{A}$ be a $T$-equivariant locally
free sheaf on $X$. Let $m_x$ be the maximal ideal of the local
ring at $x$, ${\mathcal{A}_x}$ the stalk at $x,$ and
${\mathcal{A}(x)} = {\mathcal{A}_x}/m_x{\mathcal{A}_x}$. If $x$ is
a $T$-fixed point then $T$ acts on  ${\mathcal{A}(x)}.$
%If $M$ is a $T$-module and $(r,s) \in \mathbb{Z}^2$, set $M_{r,s} = \{m \in
%M|\tau \cdot m = q^r t^s m\}$.
Suppose that $\tau$ acts on with eigenvalues $v_x(q,t)$ and
$w_x(q,t)$ on the cotangent space $m_x/m_x^2$ at $x.$ Then we have
\begin{equation} \label{ABL}
\sum_{i \geq 0}(-1)^i \mathcal{H}_{H^i(X,\mathcal{A})} (q,t) =
\sum_{x\in X^T}  \mathcal{H}_{{\mathcal{A}(x)}} (q,t)/(1 -
v_x(q,t))(1 - w_x(q,t)).
\end{equation}
\end{thm}

\noindent{\bf Proof.}
 See [AB] for the compact case.  This form
 is a special case
 of [Ha, Theorem 3.1].\\

The denominator in (\ref{ABL}) is equal to $det_{m_x/m_x^2} (1 -
\tau).$ We apply Theorem \ref{ABLT} when $\mathcal{A}$ is the
invertible sheaf corresponding to the divisor $D = D({\bf b})$  on
$X$ as in (\ref{eq5}).

\begin{lem} \label{NL1}
If ${\bf b} \in \mathbb{N}^{n-1}$ and $D = D({\bf b}),$ then
$\mathcal{O}(D)(X_{i + 1})$ is a free
   $\mathcal{O}(X_{i + 1})$-module
  generated by $x^{m_1}z^{m_2}$ where
\begin{equation}
m_1 = \sum_{j = n - i}^{n-1}(n - j) b_j , \; \quad \quad m_2 =
-\sum_{j = n - i}^{n-1}b_j.
\end{equation}
\end{lem}
\noindent {\bf Proof.} We have  $D = \sum_{i = 0}^{n}
 a_i D_i,$
 where \begin{equation} \label{eq4}
a_k = \sum_{j = n + 1 - k}^{n-1}(j +  k - n)b_j .
\end{equation}

 As in [F, page
66] we set
\[ P_D(\sigma_{i + 1}) = \{u \in M_\mathbb{R}| <u,v_j> \;\; \geq -a_j \; \;
\mbox{for} \; \; j = i, i + 1  \} .\] Then
\[H^0(X_{i + 1}, \mathcal{O}(D))
= \oplus_{u \in P_D(\sigma_{i + 1}) \cap M} \mathbb{C} x^{u_1}
z^{u_2} .\] Note that if $u = (u_1, u_2) \in M_\mathbb{R},$ then
$u \in P_D(\sigma_{i + 1})$ if and only if
\[u_1 + iu_2 \geq -a_{i } \]
and \[u_1 + (i + 1)u_2 \geq -a_{i + 1} .\] It follows that \;
$\mathcal{O}(D)(X_{i + 1})$ is freely
  generated over $\mathcal{O}(X_{i + 1})$
by $x^{m_1}z^{m_2}$ where $m_1 + im_2 = -a_{i } $ and $m_1 + (i +
1)m_2 = -a_{i + 1},$ that is

\begin{equation}
m_1 = ia_{i + 1} - (i + 1)a_{i }  \; \quad  \mbox{and } \quad m_2
= a_{i } - a_{i + 1}.
\end{equation}
Now the result follows easily from (\ref{eq4}).\\

\begin{thm}
If ${\bf b} \in \mathbb{N}^{n-1}$ and $N({\bf b}) = H^0(X,
\mathcal{O}(D({\bf b}))),$ then
\begin{equation} \label{H}
\mathcal{H}_{N({\bf b})}(q, t) = \sum_{i  = 0}^{n-1} \frac{
q^{\sum_{j = n - i}^{n-1} (n - j) b_j } \; t^{-\sum_{j = n -
i}^{n-1} jb_j}}{(1 - t^{i + 1-n}q^{i + 1})(1 - t^{n-i}q^{-i})}.
\end{equation}
\end{thm}
\noindent {\bf Proof.} By Lemma \ref{NL4}, the higher cohomology
of $\mathcal{O}(D)$ vanishes, so the left side of (\ref{ABL})
equals $\mathcal{H}_{N({\bf b})}(q, t).$ On the other hand from
(\ref{la6}) and Lemma \ref{NL1} we have
\begin{equation}
\mathcal{H}_{\mathcal{A}(p_{i + 1})} = q^{m_1}t^{m_1 + nm_2} =
%q^{-\sum_{j = i}^{n-1} jb_j } \; t^{\sum_{j = i}^{n-1}(n - j)b_j}
q^{\sum_{j = n - i}^{n-1} (n - j) b_j } \; t^{-\sum_{j = n -
i}^{n-1} jb_j},
\end{equation}
and by (\ref{m}) that the determinant of $(1 - \tau)$ on the
cotangent space at $p_{i + 1}$ is \[ (1 - t^{i + 1-n}q^{i + 1})(1
- t^{n - i}q^{-i}).\] Therefore the result follows from Theorem
\ref{ABLT}.\\

\section{Multi-homogeneous coordinate rings}
Let $\mathcal{R} = \bigoplus_{{\bf i}\in \mathbb{N}^m}
\mathcal{R}({\bf i})$ be a $\mathbb{N}^m$-graded algebra.  We
introduce several abelian categories that we will need. First we
write $\mathcal{R}\lGr$ for the category of $\mathbb{Z}^m$-graded
$\mathcal{R}$--modules with degree zero homomorphisms, and
$\mathcal{R}\lgr$ for the full subcategory of $\mathcal{R}\lGr$
consisting finitely generated graded $\mathcal{R}$--modules.

If  ${\bf i, j} \in \mathbb{Z}^m$ we write ${\bf i} \geq {\bf j}$
if ${\bf i} - {\bf j} \in \mathbb{N}^m.$ We say that a property
(P) holds for {\it large enough} ${\bf b} \in \mathbb{Z}^m$ if
there exists ${\bf j}
 \in \mathbb{Z}^m$ such that (P) holds for all ${\bf b} \geq {\bf j}$. A finitely
generated graded $\mathcal{R}$--module ${M} = \bigoplus_{{\bf
i}\in \mathbb{Z}^m} {M}({\bf i})$ is {\it bounded} if ${M}({\bf
i}) = 0$ for large enough ${\bf i}$, and a graded
$\mathcal{R}$--module is {\it torsion} if it is the direct limit
of a system of bounded modules. We denote the Serre subcategories
of $\mathcal{R}\lgr$ and $\mathcal{R}\lGr$ consisting of bounded
and torsion modules by $ \mathcal{R} \ltors$ and $ \mathcal{R}
\lTors$ respectively. Finally, we define the quotient categories
$$\mathcal{R}\mathrm{-Qgr} = \mathcal{R}\lGr / \mathcal{R}
\lTors \qquad \mathcal{R}\mathrm{-qgr} = \mathcal{R}\lgr /
\mathcal{R} \ltors.$$

Let $S_{\bf b} = H^0 (X, \mathcal{O}(D({\bf b})))$ as in Section 2
and set
 $S = \bigoplus_{{\bf b}  \in \mathbb{N}^{n-1}} S_{\bf b}.$ We
show that there is an equivalence of categories between
$S\mathrm{-qgr}$ and $\coh (X),$ the category of coherent sheaves
on $X$. To do this adapt the proof of an equivalence of categories
for
 twisted homogeneous coordinate rings, and for
 twisted multi-homogeneous coordinate rings obtained in [AV, Theorem
 1.3] and [Chan, Theorem 3.4] respectively. The proof works in greater generality.

Let $X$ be a variety, projective over $\spc S$ with $S$ a a
noetherian $k$-algebra . Assume we have a system of invertible
sheaves $(\LL_1, \ldots, \LL_m)$ on $X$ with the following
properties:
\begin{enumerate}
\item There exists $N \gg 0$ such that for all $n\geq N$ we have
$H^p (X, \LL_i^n) = 0$ for each $i$ and for all positive integers
$p$. (Note that by Grothendieck's vanishing theorem we need only
consider $p\leq \dim X$, so for each $1\leq i \leq m$ and for each
$1\leq p\leq \dim X$ we can look for a positive integer $N(i,p)$
with the cohomology vanishing property above, and then set $N$ to
be the maximum value amongst the $N(i,p)$.)

\item For a coherent sheaf $\mathcal{F}$ on X there exists a
sequence of integers $(c_1, \ldots , c_m)$ such that for all
$(b_1, \ldots , b_m) \geq (c_1, \ldots, c_m),$
%(meaning that $b_j\geq c_j$ for all $1\leq j \leq m$)
 we have $H^p(X,
\mathcal{F}\otimes \LL_1^{b_1}\cdots \LL_m^{b_m}) = 0$ for all
positive integers $p$.
\end{enumerate}
We call such a system an {\it ample system}. %(although this is
%probably not standard).
To ease notation we write
  \begin{equation} \label{Lb} \LL ({\bf
b}) = \LL^{b_1}_1 \cdots \LL^{b_m}_m
\end{equation}
 for ${\bf b} = (b_1,\ldots ,b_m)
\in \Z^m.$  Define the following graded algebras:
$$\mathcal{L} = \bigoplus_{{\bf b}\in \Z^m}
\LL({\bf b}); \quad A = \bigoplus_{{\bf b}\in \NN^m} H^0 ( X,
\LL({\bf b})).$$

Then $\mathcal{L}$ is a strongly $\Z^m$--graded
$\mathcal{O}_X$--algebra  with
 $\mathcal{L}(0) =
\mathcal{O}_X.$ The algebra $A$ is an $\NN^m$--graded
$\mathcal{O}_X(X)$--algebra.  A {\it graded} $\mathcal{L}$--module
$\mathcal{M}$ has the form
$$\mathcal{M} = \bigoplus_{{\bf i}\in \Z^m} \mathcal{M}({\bf i})$$
 with each $\mathcal{M}({\bf i})$
a quasi-coherent $\OO_X$--sheaf, such that $\mathcal{L}({\bf i})
\cdot \mathcal{M}({\bf j}) \subseteq \mathcal{M}({\bf i}+{\bf
j}).$ We have then the categories $\mathcal{L}\lGr$ of graded
$\mathcal{L}$--modules with degree zero homomorphisms, and its
subcategory $\mathcal{L}\lgr$ consisting of finitely generated
$\mathcal{L}$--modules.

\begin{thm} \label{ECQ}
There is an equivalence of categories between $\Qcoh X$ and $
A\mathrm{-Qgr}$ which restricts to an equivalence between $\coh X$
and $ A\mathrm{-qgr}$.
\end{thm}
The proof goes in two stages.  First we have
\begin{lem} \label{Dade}
The correspondences
\[ \begin{array}{ccc}
       \mathcal{M}     & \longrightarrow & \mathcal{L}  \otimes_{\mathcal{O}_X} \mathcal{M} \\
        \mathcal{N}(0) & \longleftarrow  &      \mathcal{N}
   \end{array}  \]
 induce inverse equivalences of categories between
 $\Qcoh X$ and $\mathcal{L}\lGr.$\\
\end{lem}
 \noindent {\bf Proof.}
 The proof is easily
 adapted from the corresponding result about strongly graded rings,
 [D].
%\end{proof}
\begin{lem} \label{AVlemma}
Let $\mathcal{F}$ be a coherent sheaf on $X$. Then for large
enough ${\bf i}$, $\mathcal{F}\otimes \LL({\bf i})$ is generated
by global sections.
\end{lem}
\noindent {\bf Proof.} (Adapted from [Hart, Proposition III.5.3].)
Let $P\in X$ be a closed point and $\mathcal{I}_P$ its ideal
sheaf. We have a short exact sequence
$$
0 \longrightarrow \mathcal{I}_P \mathcal{F} \longrightarrow
\mathcal{F} \longrightarrow \mathcal{F}\otimes k(P)
\longrightarrow 0.$$ Tensoring this by $\LL({\bf i})$ gives
$$
0 \longrightarrow \mathcal{I}_P \mathcal{F} \otimes \LL({\bf i})
\longrightarrow \mathcal{F} \otimes \LL({\bf i}) \longrightarrow
\mathcal{F}\otimes \LL({\bf i})\otimes k(P) \longrightarrow 0.$$
By hypothesis we can find ${\bf j}$ such that $H^1 (X,
\mathcal{I}_P\mathcal{F}\otimes \LL({\bf i})) = 0$ for all ${\bf
i}\geq {\bf j}$. Therefore the mapping $$H^0( X , \mathcal{F}
\otimes \LL({\bf i})) \longrightarrow H^0 ( X, \mathcal{F}\otimes
\LL({\bf i}) \otimes k(P))$$ is surjective for all ${\bf i}\geq
{\bf j}$. So, by Nakayama's lemma, the stalk at $P$ of
$\mathcal{F}\otimes \mathcal{L}({\bf i})$ is generated by global
sections. Hence, there exists an open neighbourhood $U$ of $P$
such that $\mathcal{F} \otimes \LL({\bf i})$ is generated by
global sections on all of $U$. However this open set depends on
the choice of ${\bf i}$. So pick the $U$ corresponding to ${\bf
j}$.

If we refine the argument in the above paragraph to work with
$\mathcal{F} = \OO_X$ we can see (using property (1) of our ample
system) that we can find neighbourhoods $V_1, \ldots , V_m$ of $P$
and integers $\ell_1, \ldots , \ell_m$ such that for each $i$,
$\LL_i^{n_i}$ is generated by global sections on $V_i$ for all
$n_i \geq \ell_i$.

Therefore on $U\cap V_1 \cap \ldots \cap V_m$ the sheaf
$\mathcal{F}\otimes \LL({\bf j}) \otimes \LL({\bf n})$ is
generated by its global sections for any $\bf n \geq l.$ Covering
$X$ with a finite number of open sets of the above form proves the
lemma.\\
%\end{proof}

To complete the proof of Theorem \ref{ECQ} we consider the
functors

\[F: A\mathrm{-mod} \longrightarrow \mathcal{L}\mathrm{-Grmod} : \quad M \longrightarrow
\mathcal{L} \otimes_{A} M\]
\[G: \mathcal{L}\mathrm{-Grmod} \longrightarrow
A\mathrm{-mod}:  \quad \mathcal{M} \longrightarrow
H^0(X,\mathcal{M})_{\geq 0} = \bigoplus_{{\bf i}\in \mathbb{N}^m}
H^0(X,\mathcal{M}({\bf i})). \]

\begin{lem} \label{NL10}
The functors $F$ and $G$ induce an equivalence of categories
between $ \mathcal{L}-\mathrm{Grmod}$ and $ A\mathrm{-Qgr}.$
\end{lem}
 \noindent {\bf Proof.} The proof in [AV] goes through verbatim using Lemma \ref{AVlemma}
 in place of [AV, 3.2 (ii)].

 \noindent {\bf Proof
of Theorem \ref{ECQ}.} This now follows by combining Lemmas
\ref{Dade} and \ref{NL10}. In particular the equivalences are
given as follows:
$$ \Qcoh X \longrightarrow A\mathrm{-Qgr} \quad : \quad \mathcal{F}
\mapsto \bigoplus_{{\bf i}\in \NN^n} H^0 (X, \mathcal{F}\otimes
\LL({\bf i})),$$ and
$$  A\mathrm{-Qgr}  \longrightarrow \Qcoh X \quad : \quad M
\longrightarrow (\mathcal{L}\otimes_A M)(0).$$
%Here the subscript$0$ denotes the degree $0$ part of a $\Z^m$--graded
%$\mathcal{L}$--module.\\

Now suppose  that $X$ is the minimal resolution the cyclic
singularity of type $A_{n-1}$ as in Section 2. The map $X
\longrightarrow \mathbb{C}^2/\Gamma$ is a projective morphism,
since it is obtained by successive blow-ups.
 For $1\leq i \leq
n-1$ set $\LL_i = \mathcal{O}(D(i)).$  Then if we define $\LL({\bf
b})$ for ${\bf b} \in \Z^{n-1}$ as in equation (\ref{Lb}),  we
have $\LL({\bf b}) = \OO(D(\bf{b}))$.

\begin{lem} \label{IS}
The invertible sheaves $\mathcal{O}(D(i))$ for  $1\leq i \leq n-1$
form an ample system on $X.$
\end{lem}

\noindent {\bf Proof.} Property (1) in the definition of an ample
system holds by Lemma \ref{NL4}.
 To check property (2) suppose that $\mathcal{F}$ be a coherent sheaf on $X$.
 We show that for
large enough ${\bf i}$ we have $H^p(X, \mathcal{F}\otimes \LL({\bf
i})) = 0$ for all positive $p$. By [Hart, Corollary II. 5.18] the
coherent sheaf $\mathcal{F}$ can be written as a quotient,
$\mathcal{E}/\mathcal{G}$, of a finite direct sum of twisted
structure sheaves. Now by [F, Proposition 3.4], the group of
$T$--Cartier divisors maps onto the Picard group of $X$, and
because $X$ is smooth every $T$--Cartier divisor is $T$--Weil.
Hence every invertible sheaf (and in particular any twisted
structure sheaf) has the form $\LL({\bf k})$ for some ${\bf k} \in
\Z^{n-1}$.

We have a short exact sequence $$0\longrightarrow
\mathcal{G}\otimes \LL({\bf i}) \longrightarrow \mathcal{E}\otimes
\LL({\bf i}) \longrightarrow \mathcal{F}\otimes \LL({\bf i})
\longrightarrow 0$$ which leads to a long exact sequence in
cohomology. For large enough ${\bf i}$ we can force the vanishing
of $H^p (X , \mathcal{E} \otimes \LL({\bf i}))$ for positive $p$
so we can deduce that $$H^p (X, \mathcal{F}\otimes \LL({\bf i}))
\cong H^{p+1}(X, \mathcal{G} \otimes \LL({\bf i})).$$ Since
$\mathcal{G}$ is coherent we have by induction that the cohomology
group on the right hand side vanishes, as required.

\section{$\mathbb{Z}^{m}$-algebras.}\label{ZA}

We need a routine generalization of the notion of a
$\mathbb{Z}$-algebra introduced in [BP]. Let $A$ be a free abelian
group
 of rank $m$.  Thus $A \cong \mathbb{Z}^{m},$ and we write
 $A^+$ for the submonoid of $A$ corresponding to $\mathbb{N}^{m}.$
  We define a partial order on $A$ by writing ${\bf i} \geq {\bf
j}$ if ${\bf i} - {\bf j} \in A^+.$
    A {\it (triangular) $\mathbb{Z}^{m}$-algebra} is an algebra
    $B = \bigoplus B({{\bf i}, {\bf j}})$
    where each $B({{\bf i}, {\bf j}})$ is an additive subgroup of $B,$ and the sum is over
    all ${\bf i}, {\bf j} \in A$
    %\mathbb{C}^{m}$
    such that ${\bf i} \geq {\bf j}.$  The multiplication in $B$
resembles matrix multiplication, that is we have $B({{\bf i}, {\bf
j}}) B({{\bf p}, {\bf q}}) = 0$ if ${\bf j} \neq {\bf p},$ and
\begin{equation} \label{cond}
B({\bf i},{\bf j}) B({\bf j},{\bf l}) \subseteq B({\bf i},{\bf
l}).
\end{equation} whenever ${\bf i} \geq {\bf j} \geq {\bf l}.$

There are two kinds of $\mathbb{Z}^m
 $-algebras that will be of interest to us.  Suppose
first that $\mathcal{R} = \bigoplus_{{\bf b}  \in \mathbb{N}^{m}}
\mathcal{R}_{\bf b}$ is an $\mathbb{N}^m$-graded algebra and set
$\widehat{\mathcal{R}} = \bigoplus_{{\bf i} \geq {\bf j} \geq 0}
\widehat{\mathcal{R}}({\bf i},{\bf j})$ where
$\widehat{\mathcal{R}}({\bf i},{\bf j}) = \mathcal{R}_{{\bf i}
-{\bf j}}.$ The multiplication in $\widehat{\mathcal{R}}$ is
induced from that in $\mathcal{R}.$ We call
$\widehat{\mathcal{R}}$ the $\mathbb{Z}^m$-{\it algebra arising
from} $\mathcal{R}.$

Let $B$ be a $\mathbb{Z}^m
 $-algebra.  We consider the category $B\mathrm{-Grmod}$ of
$\mathbb{N}^m
 $-graded left $B$-modules $M =
\bigoplus_{{\bf j} \geq 0}
 M({\bf j})$ such that
$B({\bf i},{\bf j}) M({\bf j}) \subseteq M({\bf j})$ for all ${\bf
i} \geq {\bf j}$ and $B({\bf i},{\bf j}) M({\bf k}) = 0$ if ${\bf
k} \neq {\bf j}.$ Morphisms in $B\mathrm{-Grmod}$ are graded
homomorphisms of degree zero. It is now clear how to define the
categories $B\mathrm{-grmod}$, $B\mathrm{-Qgr}$ and
$B\mathrm{-qgr}$ by analogy with the definitions we made for
$\mathbb{N}^m$-graded algebras.

Returning to the algebras $\mathcal{R}$ and
$\widehat{\mathcal{R}}$, let $\mathcal{R}\mathrm{-Grmod}_{\geq 0
}$ be the full subcategory of $\mathcal{R}\mathrm{-Grmod}$
consisting of all $\mathbb{N}^m$-graded modules.  We define the
categories $\widehat{\mathcal{R}} \mathrm{-Qgr}$ and
$\widehat{\mathcal{R}}
   \mathrm{-qgr}$ in the obvious way. For $M$ an object in
$\mathcal{R}\mathrm{-Grmod}$ we set $M_{\geq 0} = \bigoplus_{{\bf
j} \geq 0} M({\bf j}) $ and write $\pi(M)$ for the image of $M$ in
$\mathcal{R} - Qgr.$
\begin{lem} \label{ZS}
(1) The identity map yields equivalences of categories
$$\mathcal{R}\mathrm{-Grmod}_{\geq 0} \longrightarrow
\widehat{\mathcal{R}}\mathrm{-Grmod}, \;
\mathcal{R}-\mathrm{grmod}_{\geq 0} \longrightarrow
\widehat{\mathcal{R}}-\mathrm{grmod}.$$
 (2) The equivalences
in (1) induce equivalences
\[ \mathcal{R}  \mathrm{-Qgr} \longrightarrow \widehat{\mathcal{R}} \mathrm{-Qgr}
   \quad \mathcal{R} \mathrm{qgr} \longrightarrow \widehat{\mathcal{R}}
   \mathrm{-qgr}\]
\end{lem}
\noindent {\bf Proof.}  (1) is immediate. For (2)
%write $\pi(M)$ for the image in $\mathcal{R} - Qgr$ of $M \in \mathcal{R}\mathrm{-Grmod}$.  T
the only point to note is that
$\pi(M) = \pi(M_{\geq 0})$, since $M/M_{\geq 0}$ is torsion.\\

 For the second class  of  $\Z^m$-algebras,
let ${\bf w_1,\ldots ,w_{m}}$ be linearly independent elements of
$\mathbb{C}^{m},$ and set $A = \oplus_{i = 1}^{m}\mathbb{Z}{\bf
w_i}, \; A^+ = \oplus_{i = 1}^{m}\mathbb{N}{\bf w_i}.$  Suppose
that $Q$ is a $\mathbb{C}$-algebra
     %Noetherian domain with Goldie quotient     ring $Q$
     and that for ${\bf i, j} \in \mathbb{C}^{m}$
 %\oplus_{i = 1}^{m}\mathbb{N}{\bf w_i}$
 with
${\bf i} \geq {\bf j}$ we are given  subrings $R_{{\bf i}},
R_{{\bf j}}$ and $R_{{\bf i}}- R_{{\bf j}}$ sub-bimodules $B({\bf
i},{\bf j})$ of $Q$ such that $B({\bf i},{\bf i}) = R_{\bf i}$ and
(\ref{cond}) holds.
%The $\mathbb{Z}^{m}{}$-algebra we are interested in is $\oplus_{{\bf
%i}\geq {\bf j}} B({\bf i},{\bf j})$.

Fix ${\bf k} \in \mathbb{C}^{m}$ and set
  \[R({\bf k}) = \oplus_{{\bf i} \geq {\bf j} \geq {\bf k}}
                B({\bf i},{\bf j}).\]
 If the bimodules $B({\bf i},{\bf j})$ in this sum induce a Morita equivalence
 between $R_{\bf i}$ and $R_{\bf j}$, we say that $R({\bf k})$
  is {\it the
 Morita $\mathbb{Z}^{m}$-algebra associated to the data}
  $(R_{\bf i}, B({\bf i}, {\bf j}))$.  For this to happen it is
  necessary that equality holds in (\ref{cond}).

\begin{lem} \label{Mor}

 Suppose that $R({\bf k})$
is the Morita $\Z^m$-algebra associated to the data
 $R_{\bf i}, B({\bf i,j} ),$
 with $R = R_{\bf k}$ Noetherian, then
\\
 (1) each finitely generated left $R({\bf k})$-module is
Noetherian, \\
(2) the association $\phi:M \longrightarrow \oplus_{ {\bf i} \geq
{\bf k}} B({\bf i},{\bf k})\otimes_{R} M$ induces an equivalence
of categories between $R \mathrm{-mod}$ and $R({\bf
k})\mathrm{-qgr}.$
\end{lem} \noindent
  {\bf Proof.} Part (2) of the lemma follows the proof of [GS1, Lemma
5.5 (2)] verbatim, so we have only to check part (1). As in the
proof of the [GS1, Lemma 5.5] it is enough to show that
$$ M = \bigoplus_{{\bf i}\geq {\bf b}}{B({\bf i,b} )}$$ is noetherian.
So let $L = \bigoplus_{{\bf i}\geq {\bf b}}L({\bf i}) \subseteq M$
be a graded submodule. For ${\bf i}\geq {\bf b}$ set
$$X({\bf i}) = B({\bf i,b} )^* \otimes
L({\bf i}) \subseteq R_{\bf b} .$$ As $R_{{\bf b}}$ is noetherian
we have $\sum_{{\bf i}\geq {\bf b}} X({\bf i}) = \sum_{{\bf j}\geq
{\bf i}\geq {\bf b}} X({\bf i})$ for some ${\bf j}$. Observe that
there are only finitely many values of ${\bf i}$ between ${\bf b}$
and ${\bf j}$. Therefore, if ${\bf k}\geq {\bf j}$ we have
$$L({\bf k}) = B({\bf k,b} ) X({\bf k}) \subseteq
\sum_{{\bf j}\geq {\bf i}\geq {\bf b}}B({\bf k,b} ) X({\bf i})
=\sum_{{\bf j}\geq {\bf i}\geq {\bf b}}B({\bf k,i} )B({\bf i,b} )
 X({\bf i}) =\sum_{{\bf j}\geq {\bf i}\geq {\bf b}}B({\bf k,i}
)L({\bf i}).$$

Now fix $1\leq t \leq m$. For each $a$ such that $b_t \leq a <
j_t$ we have a $\Z^{m-1}$--algebra $R^{(t,a)}$ defined as
$$R^{(t,a)} = \bigoplus_{{\bf i} \geq {\bf j}, i_t=j_t
= a} B({\bf i,j} ).$$

The $R^{(t,a)}$--module
$$\bigoplus_{{\bf i} \geq {\bf b}, i_t = a} B({\bf
i,b} )$$ is finitely generated  noetherian by induction on $m$.
Therefore we can find a finite set $I(t,a) \subseteq \{{\bf i} :
i_t = a\}$ such that for all ${\bf k}$ with $k_t = a$ we have
$$L({\bf k}) \subseteq \sum_{{\bf i}\in I(t,a)}
B({\bf k, i} )L({\bf i}).$$

Thus $L$ is generated by the $L({\bf i})$ where ${\bf i}$ belongs
to the finite set $\{ {\bf i} : {\bf b}\leq {\bf i} \leq {\bf j}\}
\cup \bigcup_{1\leq t \leq m, b_t \leq a < j_t} I(t,a)$. Since
each $L({\bf i})$ is finitely generated as a $ R_{{\bf i}}$-module
the proof is finished.\\
%\end{proof}

The $\mathbb{Z}^{m}$-algebras we require can be constructed from
certain bimodules that we call {\it basic.}
%Suppose we are given only the
These are the bimodules $B({\bf j} + {\bf w_p},{\bf j})$ for ${\bf
j} \in \mathbb{C}^{m}$ and $1 \leq p \leq m.$ The basic bimodules
are required to satisfy a suitable compatibility condition. Namely
for all ${\bf j} \in \mathbb{C}^{m}$ and $1 \leq p, q \leq m$ we
require
\begin{equation} \label{ind}
B({\bf j}  + {\bf w_p + w_q}, {\bf j} + {\bf w_p}) B({\bf j} +
{\bf w_p} ,{\bf j} ) = B({\bf j}  + {\bf w_p + w_q}, {\bf j} +
{\bf w_q}) B({\bf j} + {\bf w_q} ,{\bf j} ).
\end{equation}

Then we define $B({\bf j},{\bf k})$ for ${\bf j} \geq  {\bf k}$ as
follows.  Choose ${\bf r_0}, \ldots , {\bf r_s} \in
%\oplus_{i =1}^{m}\mathbb{N}w_i
A$ such that ${\bf r_0} = {\bf k}$, ${\bf r_s} ={\bf j}$ and ${\bf
r_i} =  {\bf r_{i-1} + w_{t(i)}}$ where $t(i) \in \{1, \ldots m
\}$ for $1 \leq i \leq s.$ Then set
\begin{equation} \label{la7}
 B({\bf j} , {\bf k}) = B({\bf r_s},{\bf
r_{s-1}}) \ldots B({\bf r_i}, {\bf r_{i-1}}) \ldots B({\bf r_1} ,
{\bf r_0)}.
\end{equation}

By equation (\ref{ind}) this definition is independent of the
choice of the ${\bf r_i},$ and it is clear that (\ref{cond}) holds.\\

  \section{Morita Theory for Spherical Subalgebras.}\label{me}
 \indent Suppose that $\Gamma = (\gamma)$ is
cyclic of order $n$ and that $\Gamma$ acts on the first Weyl
algebra $\mathbb{C}[\partial, y]$ so that in the crossed product
we have
\[ y\gamma = \omega \gamma y, \quad \gamma \partial = \omega \partial \gamma\]
where $\omega = e^{2\pi i/n}.$  We do some computations in $Q =
\mathbb{C}[\partial, y^{\pm 1}]*\Gamma.$  For $0 \leq i \leq n-1$
set
\[e_i = (1/n)\sum_{j = 0}^{n-1}(\omega^i \gamma)^j\] and note that
$ ye_{i} =  e_{i+1}y$ and $ e_{i}\partial = \partial  e_{i+1},$
 where the indices are read mod $n$. Fix ${\bf k} = (k_1, \ldots, k_{n-1}) \in \mathbb{C}^{n-1}$,
 set $k_0 = 0$, and then define $k_j$ for $j \in \mathbb{Z}$ by $k_j = k_i$
 where $j \equiv i \mod n$ and $0 \leq i \leq n-1$. Then set $ e =
 e_0$ and
 \[d_{\bf k} = \partial
- y^{-1} \sum_{i = 1}^{n-1}k_i e_i.\] Let $H_{\bf k}$ be the
subalgebra of the crossed product generated by $d_{\bf k}, y$ and
$\Gamma,$ and set $U_{\bf k}= eH_{\bf k}e.$ We define
 \[ \theta = y\partial = yd_{\bf k}+ \sum_{i = 1}^{n-1} k_ie_i \in H_{\bf k}.\]% $\theta = y \partial $.
Note that  $Q = H_{\bf k}[y^{-1}]$ for all ${\bf k}.$ By induction
we have
\begin{equation} \label{Ty1}
 d_{\bf k}^p y^p = \prod_{i = 1}^{p} (\theta + i -
 \sum_{j = 0}^{n-1} k_{i+j}e_j),
\quad y^pd_{\bf k}^p  = \prod_{i = 0}^{p-1} (\theta - i -
 \sum_{j = 0}^{n-1} k_{j-i}e_j).
\end{equation}
Since the $e_i$ are orthogonal idempotents which commute with
$\theta,$ it follows that
\begin{equation} \label{Ty2}
e_j d_{\bf k}^p y^p = e_j\prod_{i = 1}^{p} (\theta + i - k_{i+j}),
\quad e_j y^pd_{\bf k}^p  = e_j \prod_{i = 0}^{p-1} (\theta - i -
  k_{j-i}).
\end{equation}

\indent Let ${\bf v_1}, \ldots {,\bf v_{n-1}}$ be the standard
basis for
    $\mathbb{C}^{n-1},$ set ${\bf w_p} = n \sum_{i = 1}^{p}{\bf
    v_i},A = \oplus_{i = 1}^{n-1}\mathbb{Z}{\bf w_i},$ and $A^+ = \oplus_{i =
1}^{n-1}\mathbb{N}{\bf w_i}.$ If ${\bf j} \geq {\bf k},$ then
${\bf j} -{\bf k} \in \sum_{j = 1}^{n-1}\mathbb{N} {\bf w_j} . $
It is convenient to define $F({\bf b}) = \sum_{j = 1}^{n-1} b_j
{\bf w_j} = n\sum_{i = 1}^{n-1} (\sum_{j = i}^{n-1} b_j) {\bf
v_j}.$ Then we have isomorphisms
\begin{equation}\label{1001} \nonumber
Pic(X) \longleftarrow  \mathbb{Z}^{n-1}
 \longrightarrow \oplus_{i =
1}^{n-1}\mathbb{Z}{\bf w_i}
\end{equation}
\begin{equation}\label{1002}
D({\bf b}) \;  \longleftarrow {\bf b} \; \longrightarrow F({\bf
b}).
\end{equation}

\begin{lem} \label{bimod}
 Fix $p \in \mathbb{N}$ with $1 \leq p \leq n-1$, and set ${\bf k'} = {\bf k} + {\bf w_p}.$
Then
\begin{equation}\label{Bimod}
 y^p e H_{\bf k}e =  e_p H_{\bf k'}e_p y^{p}.
\end{equation}
\end{lem}
\noindent
 {\bf Proof.} To simplify notation write $d = d_{\bf k}$ and $d_1 =
d_{\bf k'}.$ We have

\[ e_p \prod_{i = 1}^{p} (\theta + i -p
- k_{i})  =  e_p \prod_{j = n-p+1}^{n} (\theta + j - k'_{p+j})  \]
and
\[ e_p \prod_{i = p+1}^{n} (\theta + i -p
- k_{i})  =  e_p \prod_{j = 1}^{n-p} (\theta + j - k'_{p+j}) .\]
To see this set $j = i + n -p,$ (resp. $j = i-p$) in the left side
of the first (resp. second) equation above to obtain the right
side. Since $y\theta y^{-1} = \theta -1$ it follows from
(\ref{Ty2})
 that

\begin{equation} \label{Ty}
y^p e d^n y^{-p}  = e_p \prod_{i = 1}^{n} (\theta + i -p -
k_{i})y^{-n}  =  e_p \prod_{i = 1}^{n} (\theta + i -
k'_{p+i})y^{-n}  = e_p d_1^n.
\end{equation}
The left side of (\ref{Bimod}) equals $y^p e \mathbb{C} [y^n, yd,
d^{n}]$ and the right side equals $e_p \mathbb{C} [y^n, yd_1,
d_1^{n}]y^p.$ By (\ref{Ty}) we have
\begin{equation} \label{51}
 y^p e d^n = e_pd_1^n y^{p}
\end{equation}
and it easy to see that

\begin{equation*} \label{53}
  y^p e yd = e_p(yd_1 - \kappa)y^{p} ,
\end{equation*}
 for some $\kappa \in
 \mathbb{C}.$
 Therefore
  \[ y^p e y^{an} (yd)^b d^{cn} = e_p y^{an} (yd_1 - \kappa)^b d^{cn}_1 y^p.\]
The result follows from this.\\

We use Lemma \ref{bimod} to define the basic bimodules for our
$\mathbb{Z}^{n-1}{}-$algebras. Fix $p$ with $1 \leq p \leq n-1,$
and ${\bf k} \in \mathbb{C}^{n-1},$ and set ${\bf k'} ={\bf k} +
{\bf w_p}$ and
\[B({\bf k'},{\bf k}) = B_p({\bf k}) = eH_{\bf k'}e_p y^{p},\] By Lemma
\ref{bimod},
  \[ B_p({\bf k}) \cdot U_{\bf k} = e H_{\bf
k'}e_p\cdot y^p e H_{\bf k}e = eH_{{\bf k}'}e_pH_{\bf k'}e_p y^{p}
\subseteq B_p({\bf k})  \]
  It follows that $B_p({\bf k}) $
 is a $U_{\bf k'}-U_{{\bf k}}$-bimodule.
 Similarly $C_p({\bf k})  = y^{-p} e_p H_{\bf k'}
 e$ is a $U_{\bf k} - U_{\bf k'}$ bimodule and we have a Morita
 context
  \begin{equation}\label{MC}
     \left[ \begin{array}{cc}
         U_{\bf k'} & B_p({\bf k})\\
C_p({\bf k})  & U_{\bf k}
     \end{array} \right]
   \end{equation}

    \begin{thm} \label{BC}
 (1) $B_p({\bf k})C_p({\bf k}) = U_{{\bf k}'}$
provided $\{i - k_i  \}_{i=1}^p \cap \{ j - k_j  \}_{j= p+1}^n =
\emptyset.$\\
(2) $C_p({\bf k})B_p({\bf k})= U_{{\bf k}}$ provided $\{i - k_i
\}_{i=1}^p \cap \{ j + n - k_j  \}_{j= p+1}^n = \emptyset.$
\end{thm}
 \noindent{\bf Proof.}
\indent By (\ref{Ty2}) $B_p({\bf k})C_p({\bf k})$ contains the
elements
\begin{eqnarray*}
 e d^p_{\bf k'} e_p y^p \cdot y^{-p} e_p y^p e  =  e
 d^p_{\bf k'} y^p
   =  eg(\theta).
\end{eqnarray*}
and
  \begin{eqnarray*}
  e y^{n-p} d^{n-p}_{\bf k'} = e  h(\theta).
 \end{eqnarray*}
 where $g(\theta) = \Pi^p_{i=1} (\theta + i - n - k_i)$ and
$h(\theta) = \Pi^n_{j
       = p+1} (\theta + j - n  - k_j).$
Since $g$ and $h$ are relatively prime if and only if $\{i - k_i
\}_{i=1}^p \cap \{j - k_j  \}_{j= p+1}^n = \emptyset,$
this proves (1) and the proof of (2) is similar.\\

We can  express the conditions in Theorem \ref{BC}  in terms of a
root systems. We embed $\mathbb{C}^{n-1}$ in $\mathbb{C}^{n}$ as
$\mathbb{C}^{n-1} \times \{0\},$ and let ${\bf v}_1,\ldots, {\bf
v}_n$  be the standard basis for $\mathbb{C}^n.$   Define a
symmetric bilinear form $(\;,\;)$ on $\mathbb{C}^n$ by $({\bf
v}_i,{\bf v}_j) = \delta_{i,j}.$
 The set
\[ \Phi = \{ {\bf v}_i - {\bf v}_{j} | i \neq j \} \] forms a root system of
type $A_{n-1}.$ As a base for the root system  we choose $B =
\{\alpha_1,  \ldots ,\alpha_{n-1}\},$ where $\alpha_i = {\bf v}_i
- {\bf v}_{i+1}$ for $1 \leq i \leq n-1.$ Let $ \Phi^+$ denote the
corresponding set of positive roots. Given ${\bf a} \in
\mathbb{C}^{n-1},$ we set $\Phi_ {\bf a} = \{ \alpha \in \Phi |
({\bf a}, \alpha) \in \mathbb{Z} \}.$

Fix ${\bf k} = (k_1, \ldots, k_{n-1}) $, and let $\rho = (n-1,
\ldots, 2, 1) \in \mathbb{C}^{n-1}.$ Then
 $({\bf k} + \rho, {\bf v}_i - {\bf v}_j) = (k_i - i) - (k_j - j).$
Set $a_i = (n - i + k_i)/n$ for  $1 \leq i \leq n-1, a_n = 0$ and
${\bf a} = (a_1, \ldots, a_{n-1}).$ We have
  \begin{equation}\label{400}
{\bf k} + \rho = n{\bf a}.
 \end{equation}
Note also that $({\bf w_p}, \alpha_i) = n\delta_{i,p}$ for $1 \leq
i,p \leq n-1.$ Let ${\bf 1} = \sum_{i=1}^{n}{\bf v}_i.$ It follows
that $({\bf w_p} - p{\bf 1})/n$  is the $p^{th}$ fundamental
weight corresponding to the basis $B.$
 We say that ${\bf k}$ is {\it dominant} if $({\bf k} + \rho,\alpha)
> 0$ for all $\alpha \in \Phi_{\bf a} \cap \Phi^+.$
For similar definitions in the enveloping algebra context see [J,
Section 2.5].

\begin{thm} \label{ad} If
${\bf k}$ is  dominant  and ${\bf k'} = {\bf k} + {\bf w_p},$
then ${\bf k'}$ is dominant and $U_{\bf k}$ and $U_{\bf k'}$ are
Morita equivalent.
\end{thm}
 \noindent
{\bf Proof.} The Morita equivalence follows from Theorem \ref{BC},
and it is easy to check that ${\bf k'}$ is dominant.\\

%We construct a Morita $\Lambda{}-$algebra with $R_0 = U_{\bf k}.$

 \indent We check that the basic bimodules $B({\bf k'},{\bf
k}) $ satisfy (\ref{ind}).

\begin{lem} \label{NL8}
If ${\bf k'} = {\bf k} + {\bf w_p}, \; {\bf k^{''}} = {\bf k} +
{\bf w_q},$ and ${\bf k^{'''}} = {\bf k} + {\bf w_p}  + {\bf
w_q},$ then
\[ B({\bf k}^{'''},{\bf k}^{'})  B({\bf k}^{'},{\bf k}) =
B({\bf k}^{'''},{\bf k}^{''})  B({\bf k}^{''},{\bf k}).
\]
\end{lem}
 \noindent
{\bf Proof.} We first show that
\begin{equation} \label{eq1}
 B({\bf k}^{'},{\bf k}) = U_{\bf k^{'}} ed_{\bf k^{'}}^py^p  +
U_{\bf k^{'}}ey^n.
\end{equation}
 Since $e_{0}H_{\bf k^{'}}e_p = e
\mathbb{C}[d_{\bf k^{'}} , y]e_p,$ and \[ ed_{\bf k^{'}}^a y^be_p
=
 e e_{p+b - a}d_{\bf k^{'}}^a y^b,\] it follows that $ eH_{\bf k^{'}}e_p$ is spanned
by all elements of the form $ed_{\bf k^{'}}^a y^be_p$ with $a
\equiv p+b \; mod \;n.$  Equation  (\ref{eq1}) follows easily from
this. Thus
\[ B({\bf k}^{'''},{\bf k}^{'}) B({\bf k}^{'},{\bf k})  =
B({\bf k}^{'''},{\bf k}^{'})d_{\bf k^{'}}^p y^pe + B({\bf
k}^{'''},{\bf k}^{'})y^ne\]
%6 = {\bf k}) U_{\bf k^{'''}\]
Using the analog of (\ref{eq1}) for $B({\bf k}^{'''},{\bf
k}^{'}),$  we see that $ B({\bf k}^{'''},{\bf k}^{'}) B({\bf
k}^{'},{\bf k})$ is generated as a left $U_{\bf k^{'''}}$- module
by the elements \[ ed_{\bf k^{'''}}^qy^q \cdot  d_{\bf k^{'}}^p
y^p, \quad ey^n\cdot d_{\bf k^{'}}^p y^p,\]
\[ ed_{\bf k^{'''}}^qy^q \cdot  y^n,
\quad ey^n\cdot y^n.\] Similarly $B({\bf k}^{'''},{\bf k}^{''})
B({\bf k}^{''},{\bf k})$ is generated as a left $U_{\bf k^{'''}}$-
module by the elements \[ ed_{\bf k^{'''}}^py^p \cdot d_{\bf
k^{''}}^q y^q, \quad ey^n\cdot d_{\bf k^{''}}^q y^q,\]
\[ ed_{\bf k^{'''}}^py^p \cdot  y^n,
\quad ey^n\cdot y^n.\] Assume that $p \leq q.$ We have the
following identities,\\
%{\bf I need to finish this.}
\begin{eqnarray}
 e y^n d^p_{\bf k'} y^p & = & e d^p_{\bf k'''} y^{p + n}
\end{eqnarray}
and
\begin{equation}\label{Z}
 e d^q_{\bf k'''} y^q d^p_{\bf k'} y^p = e d^p_{\bf k'''} y^p d^q_{\bf k''}
 y^q.
\end{equation}
These identities follows easily from (\ref {Ty2}).  For example
both sides of (\ref {Z}) are equal to
 \[ e \Pi^p_{i=1}(\theta + i - 2n - k_i) \Pi^q_{i=1}(\theta
 + i  - n - k_i) .\]

It remains to show that
\begin{eqnarray} \label{Z1}
e d^q_{\bf k'''} y^{q+n} \in B({\bf k}^{'''},{\bf k}^{''})  B({\bf
k}^{''},{\bf k})
\end{eqnarray}
and
\begin{eqnarray} \label{Z2}
 e y^n d^q_{\bf k''} y^q \in B({\bf k}^{'''},{\bf k}^{'})  B({\bf k}^{'},{\bf
 k}).
\end{eqnarray}
By (\ref{Ty2}) $ e y^n d^q_{\bf k''}y^q$ is a
$\mathbb{C}[\theta]$-multiple of $ e y^nd^p_{\bf k'} y^p$.  This
gives (\ref{Z2}) and the
proof of (\ref{Z1}) is similar.\\

By Lemma \ref{NL8} it now makes sense to define the bimodules
$B({\bf j},{\bf k})$ for ${\bf j} \geq  {\bf k}$ using equation
(\ref{la7}).
\\
\begin{thm} \label{MLA} If
${\bf k}$ is  dominant then
  $R({\bf k})  =
\oplus_{{\bf i} \geq {\bf j} \geq {\bf k}}
                B({\bf i},{\bf j})  $
  is a Morita $\mathbb{Z}^{n-1}$-algebra
$R_{\bf k} = U_{\bf k}$.
\end{thm}

  \section{Proof of the Main Theorem.}\label{sc}
Suppose that  $R = R({\bf k})$ is as in Theorem \ref{MLA}. The
algebra  $D =  \mathbb{C}[y^{\pm1},
\partial]\ast\Gamma$ is equal to $D =  \mathbb{C}[y^{\pm1},
d_{{\bf k}}]\ast\Gamma$ for all   ${\bf k} \in \C^{n-1}.$  We
consider the differential operator filtration on $D$ defined by
$D_N = \bigoplus_{i=0}^N
\mathbb{C}[y^{\pm1}]\ast\Gamma\partial^i.$ We have $D_N =
\bigoplus_{i=0}^N \mathbb{C}[y^{\pm1}]\ast\Gamma d_{{\bf k}}^i.$
If $a =
 \sum_{j = 0}^{N}f_j(y)\partial^j \in \mathbb{C}[y^{\pm1},
\partial]$ with $f_i(y) \in \mathbb{C}[y^{\pm1}], f_N \neq 0$ we set $gr
(a) = f_N(u)v^N .$ We extend $gr$ to a linear map
 from $D =  \mathbb{C}[y^{\pm1},
\partial]\ast\Gamma$ to $\mathbb{C}[u^{\pm1},
v]\ast\Gamma$ such that $gr(\gamma a)= \gamma gr(a).$ Note in
particular that $u = gr \; y, v = gr \; d _{{\bf k}} $ and $gr \;
d_{{\bf k}}= v^n = z.$ Also

\begin{eqnarray} \label{gr}
gr(bd_{\bf k}^{n}) = gr(b)z
\end{eqnarray}
 for all
 $b \in D.$ Since $H_{\bf k}$ and $ U_{\bf k}$ are subalgebras of $D$ they have an induced
filtration, and we have
\[ gr H_{{\bf k}} \cong \mathbb{C}[u, v]\ast\Gamma, \quad gr
U_{{\bf k}} \cong   \mathbb{C}[u, v]^{\Gamma}.
\]

Similarly there is a differential operator filtration $\{ B_n({\bf
r'},{\bf r}) \}_{n \geq 0}$
 on $B({\bf r'},{\bf r})$ and we obtain a filtration on
  $R$ by setting
  $R_n  =\oplus_{{\bf r'} \geq {\bf r} \geq {\bf k}}
  B_n({\bf r'},{\bf r}).$\\

Recall the isomorphisms from equation (\ref{1002}). The key
remaining step in the proof of the main theorem is the following.

\begin{thm} \label{KRS}
Assume that  ${\bf k} \in \C^{n-1}$ is dominant  and
   ${\bf r}' \geq {\bf r} \geq {\bf k}$ with ${\bf r}' = {\bf r}+ F({\bf b}).$
Then
\[gr B({\bf r'},{\bf r})  = H^0(X,
\mathcal{O}(D({\bf b}))).\]
\end{thm}
\noindent {\bf Proof of the Main Theorem.} \\
(1)  The equivalence of categories
   $R-\mathrm{qgr} \simeq U_{\bf k}$-mod follows from Theorem \ref{MLA} and Lemma
   \ref{Mor}.\\
(2) Let $S = \bigoplus_{{\bf b}  \in \mathbb{N}^{n-1}}
H^0(X,\mathcal{O}(D({\bf b}))).$  By Theorem \ref{ECQ} and Lemma
\ref{IS} we know the category $S-\mathrm{qgr}$ is equivalent to
   $Coh(Hilb_\Gamma
\mathbb{C}^2).$  On the other hand by Theorem \ref{KRS} the
associated graded ring of the $\mathbb{Z}^{n-1}$-algebra
  $R({\bf k})$ is
\[\bigoplus_{{\bf r'} \geq {\bf r} \geq {\bf k}}H^0(X,
\mathcal{O}(D(F^{-1}({\bf r'} - {\bf r}))),\]
  and this is the  $\mathbb{Z}^{n-1}$-algebra arising
  from the $\mathbb{N}^{n-1}$-graded algebra
$S.$  Hence the result follows from Lemma \ref{ZS}.\\

It is easy to check Theorem \ref{KRS} for basic bimodules.

\begin{lem} \label{NL2}
If $1 \leq p \leq n-1,$ and ${\bf r'} = {\bf r} + {\bf w_p},$ we
have
\[ gr(B({\bf r'},{\bf r})) \cong H^0(X, \mathcal{O}(D({\bf v_p}))
\]
\end{lem}
\noindent {\bf Proof.} By (\ref{eq1}) $gr B({\bf r'},{\bf r})
\cong C^{\Gamma} u^pv^p + C^{\Gamma}u^n = C^{\Gamma}x^p +
C^{\Gamma}x^n z^{-1}.$ On the other hand Lemma \ref{NL1} implies
that
\[ H^0(X, \mathcal{O}(D({\bf v_p}))
 = C^{\Gamma} + C^{\Gamma}x^{n-p} z^{-1} .\] The isomorphism is
 multiplication by $x^{-p}.$

\begin{lem} \label{inc} For ${\bf b} \in \mathbb{N}^{n-1},$ set
$f({\bf b}) = \sum_{j = 1}^{n-1}j b_j.$ There is an injective
linear map
\begin{equation}\label{inceq}
 H^0(X, \mathcal{O}(D({\bf b})))
\longrightarrow gr \; B({\bf r} + F({\bf b}),{\bf r} )
 \end{equation}
 given by multiplication by $x^{f( {\bf b})}.$
\end{lem}
\noindent {\bf Proof.} Using Corollary \ref{NL3}, induction and
then [GS1, Lemma 6.7 (1)]
we get\\

$H^0(X, \mathcal{O}(D({\bf b} +{\bf c})))x^{f( {\bf b} +  {\bf
c})} = H^0(X, \mathcal{O}(D({\bf b})))x^{f( {\bf c})}
H^0(X,\mathcal{O}(D({\bf c})))x^{f( {\bf b})} \subseteq $\\

$gr \; B({\bf r} + F({\bf b}+ {\bf c}),{\bf r} +F({\bf c})
)gr\;B({\bf r} + F({\bf c}),{\bf r} ) \subseteq gr \;B({\bf
r}+F({\bf b} + {\bf c}),{\bf r}).$\\

 From now on we assume that  ${\bf k} \in \C^{n-1}$ is dominant  and fix
   ${\bf r'} \geq {\bf r} \geq {\bf k}$ with ${\bf r'} = {\bf r}+ F({\bf b}).$
   Note that  ${\bf r'}$ and  ${\bf r}$ are also dominant. We
   write
 $B({ \bf b})$ and $d$ instead of $B({\bf r'},{\bf r})$ and
$d_{\bf r}$  respectively
%When ${\bf r}' = {\bf r + w_p}$ we often write $B_p({\bf r})$ instead of $B({\bf r}',{\bf r})$.
We have
 \begin{equation}\label{kk}
  r_i' = r_i + n\sum_{j = i}^{n-1}b_j
 \end{equation}
 for $1 \leq i \leq n-1.$\\

We explain how the proof of Theorem \ref{KRS} reduces to a
Poincar\'{e} series computation.
%By we have an injective
%maps of $\mathbb{C}[z]$-modules \[H^0(X, \mathcal{O}(D({\bf
%b})))Q^{\bf b} \longrightarrow gr \; B({\bf b}).\]
Note that the torus $T$ acts on both the domain and target of the
map in equation (\ref{inceq}). Set $P = gr \; B({\bf b})$ and
consider the weight space decomposition $P = \bigoplus_{\chi \in
\mathbb{X}(T)} P({\chi}).$ Since $gr \; B({\bf b}) \subseteq
\mathbb{C}[u,v]$ it follows that \[\{\chi \in \mathbb{X}(T)|
P({\chi}) \neq 0\} \subseteq \mathbb{N}\chi_1 +
\mathbb{N}\chi_2.\] Because $z \in C[{n\chi_2}],$ this implies
that  $\cap_{m \geq 0}Pz^m = 0.$ Therefore by Nakayama's lemma if
$\mathcal{B}$ is a subspace of $P$ whose image in $\overline{P} =
P\otimes_{\mathbb{C}[z]} \mathbb{C}$ is a basis, then
$\mathcal{B}$ generates $P$ as a $\mathbb{C}[z]$-module.  We will
show that  $dim \overline{P}({\chi})< \infty$ for all $\chi \in
\mathbb{X}(T)$. Then since $P$ is a torsion-free
$\mathbb{C}[z]$-module, we have
 \[ \mathcal{H}_{\overline{P}}(q,t) =
 (1-t^n)\mathcal{H}_{P}(q,t). \]
 Now we view $\mathbb{C}^\times$ as the subtorus $\{(\lambda,
 \lambda^{-1})|\lambda \in \mathbb{C}^\times \}$ of $T.$  Then $\mathbb{X}(\mathbb{C}^\times) =
 \mathbb{Z}\eta,$ where $\eta =
\chi_{_1|_{_{\mathbb{C}^\times}}} =
-\chi_{_2|_{_{\mathbb{C}^\times}}} .$ For a
$\mathbb{C}^\times$-module $V$ we now have the weight space
decomposition $V = \bigoplus_{m \in \mathbb{Z}} V[{m}],$ where
\[V[{m}] = \{ v \in V| \tau \cdot v = \eta(\tau)^m v \; \mbox{for
all} \; \tau \in \mathbb{C}^\times\}.\]  If  $dim V[{m}]< \infty$
for all $m,$ we define the one variable Poincar\'{e} series

\[ \mathbb{H}_V(s) = \sum_{m \in \mathbb{Z}} \dim
V[{m}]\;s^m.\] Clearly $\overline{P}[m] = \bigoplus
\overline{P}({\chi})$ where the sum is over all  $\chi \in
\mathbb{X}(T)$ such that $\chi|_{_{\mathbb{C}^\times}} = m \eta .$
If we can show that $dim \overline{P}[m] < \infty$ for all $m \in
\mathbb{Z},$ then $dim \overline{P}({\chi})< \infty$ for all $\chi
\in \mathbb{X}(T)$ and we have
 \[ \mathcal{H}_{\overline{P}}(s,s^{-1}) =
\mathbb{H}_{\overline{P}}(s). \] Similarly if $\overline{N} =
N({\bf b})\otimes_{\mathbb{C}[z]}\mathbb{C},$ then
 \[ \mathbb{H}_{\overline{N}}(s) =
\mathcal{H}_{\overline{N}}(s,s^{-1}) = (1 -s^{-n})
\mathcal{H}_{{N}}(s,s^{-1}).
\] Therefore by equation
(\ref{H})
\begin{equation} \label{obar}
\mathbb{H}_{\overline{N}}(s) = \frac{\sum_{i = 0}^{n-1} s^{n
\sum_{j = n-i}^{n-1}b_j}}{(1 - s^{n})}.
\end{equation}

On the other hand the map in equation (\ref{inceq}) is equivariant
for the action of $\mathbb{C}^\times,$ since $x = uv \in C[{0}].$
  Hence  to prove Theorem \ref{KRS} it suffices to show that
the Poincar\'{e} series $\mathbb{H}_{\overline{P}}(s)$ is given by
equation  (\ref{obar}). We do this by developing some
representation theory of the algebras $H_{\bf r}$ and $U_{\bf r},$
and the bimodules $B({\bf
b})$.\\

The algebra $H_{\bf r}$ has a $\mathbb{Z}$-grading in which the
degrees of the
 generators satisfy
 \[ \deg_{} y = 1, \quad \deg_{}  d  = -1, \quad \deg_{}  \gamma = 0 .\]

 \indent If $M = \oplus_{\alpha \in
\mathbb{Z}} M_\alpha$ is a graded $H_{\bf r}$-module with $\dim
M_\alpha < \infty$ for all $\alpha,$  we define the Poincar\'{e}
series of $M$ to be
\[ p(M,s) = \sum_{\alpha \in \mathbb{Z}}(\dim M_\alpha) s^\alpha .
\]
\begin{lem} \label{LN}

\[ \mathbb{H}_{\overline{ gr \;B({\bf b})}}(s) =
p(B({\bf b}) \otimes_{\mathbb{C}[d^n]}\mathbb{C},s). \]
\end{lem}
\noindent {\bf Proof.} Note that  $u = gr \; y \in
\mathbb{C}[u,v][1],$ and $v = gr \; d \in \mathbb{C}[u,v][-1].$
Since $\mathbb{Z} \cong \mathbb{Z}\eta$ is the character group of
$\mathbb{C}^{\times},$ this gives an action of $\mathbb{C}^\times$
on $\mathbb{C}[y^{\pm1},
\partial]\ast\Gamma$ such that the map $gr $ is
$\mathbb{C}^\times$-equivariant.  Therefore as
$\mathbb{C}^\times$-modules
\begin{equation} \label{ne7}
 B({\bf b})  \cong gr \; B({\bf b})
.\end{equation}
 Furthermore by equation (\ref{gr})
\begin{equation} \label{ne8}
 B({\bf b})d^n \cong gr \; B({\bf b})z.
\end{equation}
The result follows since $B({\bf b})
\otimes_{\mathbb{C}[d^n]}\mathbb{C}$  is the quotient of the left
side of (\ref{ne7}) by the left side of (\ref{ne8}), and
$\overline{ gr \;B({\bf b})} = (gr \; B({\bf b}))
\otimes_{\mathbb{C}[z]}\mathbb{C}$ is the corresponding quotient of the right sides.\\

 Let
$\mathcal{O}_{\bf r}$ denote the category whose objects are
finitely generated $H_{\bf r}$-modules on which the action of $d$
is locally nilpotent.
 As in [GS1, Section 6.12] and [GGOR,
Section 2.4] we use a graded version $\mathcal{\widetilde{O}}_{\bf
r}$ of the category $\mathcal{O}_{\bf r}$.
 Objects in $\mathcal{\widetilde{O}}_{\bf r}$ are
 $\mathbb{Z}$-graded $H_{\bf r}$-modules
 which are also objects in $\mathcal{{O}}_{\bf r}$.
Morphisms in $\mathcal{O}_{\bf r}$ (resp.
$\mathcal{\widetilde{O}}_{\bf r}$) are $H_{\bf r}$-module
homomorphisms (resp. $H_{\bf r}$-module homomorphisms which are
homogeneous of degree zero).  We write Let $\mathcal{O}_{\bf r}^U$
for  the category of finitely generated $U_{\bf r}$-modules on
which the action of $d^n$ is locally nilpotent, and let
$\widetilde{\mathcal{O}}_{\bf r}^U$ denote the corresponding
category of
graded $U_{\bf r}$-modules.\\
 \indent If $M = \oplus_{\alpha \in
\mathbb{Z}} M_\alpha$ is a module in $\mathcal{\widetilde{O}}_{\bf
r}$ it follows from the local nilpotence of $d$ and finite
generation that $\dim M_\alpha < \infty$ for all $\alpha$, so
$p(M,s)$ is defined. For $\beta \in \mathbb{Z}$, the shift functor
$[\beta]$ in $\widetilde{\mathcal{O}}_{\bf r}$ is defined by
$(M[\beta])_\alpha = M_{\alpha - \beta}$.  We have
\[ p(M[\beta],s) = s^\beta p(M,s) . \]
The algebra $H_{\bf r}$ has a triangular decomposition
\[ H_{\bf r} = \mathbb{C}[y] \otimes \mathbb{C}[\Gamma] \otimes
\mathbb{C}[d]. \]

For $i = 0, \ldots, n - 1$ let $\mathbb{C}\varepsilon_i$ denote
the one-dimensional $\mathbb{C}\Gamma$-module on which $e_i$ acts
as the identity and make $\mathbb{C}\varepsilon_i$ into a
$\mathbb{C}[d]*\Gamma$-module with $d \varepsilon_i = 0.$ Then we
define the {\it standard module} $M_{\bf r}(\varepsilon_i)$ by
\[ M_{\bf r}(\varepsilon_i) = H_{\bf r}
\otimes_{\mathbb{C}[d]*\Gamma} \mathbb{C}\varepsilon_i. \] Note
that $\theta$ acts on the subspace $1 \otimes
\mathbb{C}\varepsilon_i$ of $ M_{\bf r}(\varepsilon_i)$ as
multiplication by the scalar $r_i.$
 We also define the {\it graded standard module}
 $\widetilde{M}_{\bf r}(\varepsilon_i)$ to be an isomorphic copy of
 $M_{\bf r}(\varepsilon_i)$ as an $H_{\bf r}$-module, with grading
 defined by $\deg(1 \otimes \varepsilon_i) = 0$.

 \indent To prove the main theorem we compute the Poincar\'{e}
 series of $B({\bf b})\otimes_{\mathbb{C}[d^n]}\mathbb{C}$ and show
 that it equals (\ref{obar}).  To do this we observe that
\[ B({\bf b}) \otimes_{\mathbb{C}[d^n]}\mathbb{C} \cong B({\bf b}) \otimes_{eH_{\bf
r}e} eH_{\bf r}e \otimes_{\mathbb{C}[d^n]} \mathbb{C}. \]
 So we begin our analysis with the left $H_{\bf r}$-module
 $G = H_{\bf r} e \otimes _{\mathbb{C}[d^n]} \mathbb{C}.$
This module inherits a grading from $H_{\bf r}$, in which $\deg(e
\otimes 1) = 0$.
\begin{lem}
In the Grothendieck group of the category
$\mathcal{\widetilde{O}}_{\bf r}$, we have
\begin{equation}\label{G}
[G] = [ \widetilde{M}_{\bf r}(\epsilon_0)] +
\sum^{n-1}_{i=1}[\widetilde{M}_{\bf r}(\varepsilon_{n-i})[-i]].
\end{equation}
\end{lem}
\noindent {\bf Proof.}  For $0 \leq i \leq n - 1$ set $u_i = d^i e
\otimes 1 \in G$, $N_i = H_{\bf r} u_i$ and  $N_n = 0.$  It is
easy to see that $G = \sum^{n-1}_{i=0}\mathbb{C}[y]u_i$, a free
$\mathbb{C}[y]$-module of rank $n$.    Thus $N_i/N_{i+1} \cong
\widetilde{M}_{\bf r}(\varepsilon_{n-i})[-i]$ for $0 \leq i \leq
n-1$, and the result follows from this.

\begin{klem} \label{bten} Assume that ${\bf r}$ is
 dominant.
Then as objects of the category $\mathcal{\widetilde{O}}_{\bf
r'}^U$ we have, for $1 \leq i \leq n - 1$
        \begin{equation}\label{la1}
 B({\bf b}) \otimes e \widetilde{M}_{\bf r}(\varepsilon_i) \cong
 e \widetilde{M}_{\bf r'}(\varepsilon_i)[r'_i - r_i ].
        \end{equation}
 \end{klem}

To prove this we need several preliminary results. First we state
an easy characterization of the $U_{\bf r^{'}}$-modules $e{M}_{\bf
r'}(\varepsilon_i).$
\begin{lem} \label{chmod}
Suppose that $M = \mathbb{C}[y^n]v$  is a left $U_{\bf
r^{'}}$-module which is generated by an element $v$ such that
$d_{\bf r'}^nv = 0, \; \theta v = (n - i  + r'_i)v,$ and $M$ is a
free $\mathbb{C}[y^n]$-module, then $M \cong
 e {M}_{\bf r'}(\varepsilon_i).$
 \end{lem}

In the next lemma we assume that $\mathcal{R}$ and $\mathcal{S}$
are subrings of a $\mathbb{C}$-algebra $\mathcal{Q}$ and that
$\mathcal{B}$ is an $\mathcal{R}- \mathcal{S}$ bimodule and
$\mathcal{C}$ an $\mathcal{S}-\mathcal{R}$-bimodule such that the
functors

\[\mathcal{B} \otimes \underline{\;\;\;}:\mathcal{S}\mathrm{-mod}\longrightarrow \mathcal{R}\mathrm{-mod}\]
and
\[\mathcal{C} \otimes \underline{\;\;\;}:\mathcal{R}\mathrm{-mod} \longrightarrow \mathcal{S}\mathrm{-mod}
\]
are inverse equivalences of categories.  Suppose that
$\mathcal{T}$ is a multiplicatively closed subset of both
$\mathcal{R}$ and $\mathcal{S}$ and  that $\mathcal{T}$ is an Ore
set in $\mathcal{R}.$ Assume also that $\mathcal{C}$ satisfies an
Ore condition with respect to $\mathcal{T}:$ given $t \in
\mathcal{T}$ and $c \in \mathcal{C}$ there exist $t' \in
\mathcal{T}$ and $c' \in \mathcal{C}$  such that
\begin{equation} \label{Ore}
 c't = t'c.
\end{equation}
  In the following ``torsion'' means
torsion with respect to $\mathcal{T}$.
\begin{lem} \label{easy} Let $M$ be an $\mathcal{S}$-module, and N an $\mathcal{R}$-module.\\
(1) If N is torsion so is $\mathcal{C} \otimes N.$ \\
(2) If $M$ is torsion-free so is $\mathcal{B} \otimes M.$
 \end{lem}
\noindent {\bf Proof.}  (1)  Suppose that $n \in N, t \in
\mathcal{T}$  and $t n  = 0.$ If  $c \in \mathcal{C},$ we find
$c'$ and $t'$ as in (\ref{Ore}). Then $t'(c \otimes
n) =  c' \otimes t n  = 0.$\\
(2)  Since $\mathcal{T}$ is Ore in $\mathcal{R}$ the set $N$ of
torsion elements of  $\mathcal{B} \otimes M$ forms a submodule. By
(1), $\mathcal{C} \otimes N$ is a torsion submodule of
$\mathcal{C} \otimes \mathcal{B} \otimes M,$ but by assumption
$\mathcal{C} \otimes \mathcal{B} \otimes M \cong M$ which is
torsion free.
\begin{cor}  \label{tfree}
$B_p({\bf r}) \otimes e {M}_{\bf r}(\varepsilon_i)$ is a
torsion-free
 $\mathbb{C}[y^n]$-module.
\end{cor}
\noindent {\bf Proof.}
   We apply Lemma \ref{easy} with
 $\mathcal{R} = U_{\bf r'}, \mathcal{S} = U_{\bf r}, \mathcal{B} =
 B_p({\bf r})$ and $\mathcal{T} =
 \mathbb{C}[y^n] \setminus \{0\} .$
The Ore conditions follow from the existence of the localizations.
Indeed, it is easy to see that
\begin{eqnarray*}
H_{\bf r} \otimes_{\mathbb{C}[y^n]}\mathbb{C}(y^n) & = &
\mathbb{C}[y^{\pm 1}, \partial] * \Gamma\\
U_r \otimes_{\mathbb{C}[y^n]} \mathbb{C}(y^n) & = & e
\mathbb{C}[y^{\pm 1}, \partial]* \Gamma e
\end{eqnarray*}
and
\begin{eqnarray*}
B_p({\bf r}) \otimes_{\mathbb{C}[y^n]} \mathbb{C}(y^n) & =
& e \mathbb{C}[y^{\pm 1}, \partial] * \Gamma e_p y^p \\
 & = & \mathbb{C}(y^n) \otimes_{\mathbb{C}[y^n]} B_p({\bf r}).
\end{eqnarray*}
\begin{lem} \label{NL5} Assume that ${\bf r}$
is dominant and set $B = B_p({\bf r})$.  Consider the elements
 $v = ey^n \otimes y^{n-i}\varepsilon_{i}$  and $w = e d_{\bf
r'}^p y^p \otimes y^{n-i}\varepsilon_{i}$ of $B \otimes e {M}_{\bf
r}(\varepsilon_i).$    Then\\
(1) $d_{\bf r'}^nw = 0.$\\
(2) If  $1 \leq i \leq p,$ then $d_{\bf r'}^n v = w = 0.$\\
(3) If $p + 1 \leq i \leq n-1,$ then  $y^n w = \kappa v$ for some
nonzero $\kappa \in \mathbb{C}.$\\
(4)  $B \otimes e {M}_{\bf r}(\varepsilon_i) = U_{\bf
r^{'}}v + U_{\bf r^{'}}w.$\\
\end{lem}
\noindent {\bf Proof.} We use the following identity repeatedly

\begin{equation} \label{52}
  (\theta + j - r_{j}')y^{n-i}\varepsilon_{i} = ( r_i + n -i+j -  r_{j}')y^{n-i}\varepsilon_i.
\end{equation}
Equation (\ref{52}) holds because both sides are equal to
$y^{n-i}(\theta + n -i+j -  r_{j}')\varepsilon_i.$ \\
(1) Multiplying
 (\ref{51}) on the left by $d_{\bf r'}^p$ gives $ ed_{\bf r'}^n\cdot
 d_{\bf r'}^p y^p  = ed_{\bf r'}^p y^p \cdot d_{\bf r}^n.$  Therefore
 \[ d_{\bf r'}^nw = ed_{\bf r'}^n\cdot
 d_{\bf r'}^p y^p \otimes y^{n-i}\varepsilon_{i}\]
\[=  ed_{\bf r'}^p y^p \otimes d_{\bf r}^ny^{n-i}\varepsilon_{i} =
0.\]

\noindent  (2) We first show that $y^nw = 0.$ Indeed by
(\ref{Ty2}) and (\ref{52}),

\[y^nw = ey^n \cdot d_{\bf r'}^p y^p \otimes y^{n-i}\varepsilon_{i}  \]

\[= ey^n \otimes \prod_{j = 1}^{p} (\theta + j - r_{j}')  y^{n-i}\varepsilon_{i}\]

 \[= ey^n \otimes y^{n-i}\prod_{j = 1}^{p}
 ( r_i + n -i+j -  r_{j}')\varepsilon_i =  0,\]
since the term in the product with $j = i$ is zero.\\

 It follows from Corollary \ref{tfree} that $w = 0.$  Now by
equations (\ref{Ty2}) and (\ref{52}),
\[ d_{\bf r'}^nv = ed_{\bf r'}^n y^n \otimes y^{n-i}\varepsilon_{i} =
 e\prod_{j = 1}^{n} (\theta + j - r_{j}') \otimes y^{n-i} \varepsilon_{i}\]

 \[ = ed_{\bf r'}^p y^p \otimes
\prod_{j = p+1}^{n} (\theta + j - r_{j}')y^{n-i}\varepsilon_{i} \]

  \[= ed_{\bf r'}^p y^p \otimes y^{n-i}\prod_{j = p+1}^{n}
 ( r_i + n -i+j -  r_{j}')\varepsilon_i. \]
The result follows since this is a multiple of $w.$

\noindent  (3) As in the proof of (2),

\[y^nw = ey^n \otimes y^{n-i}\prod_{j = 1}^{p}
 ( r_i + n -i+j -  r_{j}')\varepsilon_i,\]
 but now the terms $ ( r_i + n -i+j -  r_{j}')$ are
nonzero, since $i \geq p+1 > j$ and ${\bf r}$
is dominant.\\
\noindent (4)  Since $B \otimes e \widetilde{M}_{\bf
r}(\varepsilon_i) = B \otimes y^{n-i}
\mathbb{C}[y^n]\varepsilon_{i} = B \otimes y^{n-i}
\varepsilon_{i},$ the result follows from equation (\ref{eq1}).\\

\noindent {\bf Proof of the Key Lemma.} We can assume that ${\bf
r'} = {\bf r} + {\bf w_p}.$  Set $B = B_p({\bf r}).$ Suppose that
$p + 1 \leq i \leq n-1.$ By (1), (3) and (4) in Lemma \ref{NL5}
 $B \otimes e \widetilde{M}_{\bf r}(\varepsilon_i) = U_{\bf r^{'}}w,$
and $d_{\bf r'}^nw = 0.$  Hence $B \otimes e {M}_{\bf
r}(\varepsilon_i)$ =
 $\mathbb{C}[y^n]w$ and this is a free $\mathbb{C}[y^n]$-module, since it is
 torsion free by Corollary \ref{tfree}.
 It is easy to check that $\theta w = (n -i  +  r_i)w$ and that $w$ has the same
degree as $y^{n-i}\varepsilon_{i}.$ So the result in this case
 follows from Lemma \ref{chmod}.
 The proof for the case where $1 \leq i \leq p,$ is
similar using $v$ instead of $w.$\\
\\
\noindent{\bf Proof of Theorem \ref{KRS}.} By Lemma \ref{inc}
there is an inclusion $H^0(X, \mathcal{O}(D({\bf b}))) \subseteq
gr B({\bf b}) $. To show the reverse inclusion, it suffices to
show that
 \[ gr B({\bf b}) \otimes_{\mathbb{C}[z]} \mathbb{C} \]
       and
 \[ H^o(X,\mathcal{O}(D({\bf b}))) \otimes_{\mathbb{C}[z]}
    \mathbb{C} \]
  have the same Poincar\'{e} series.\\

By equations (\ref{G}) and (\ref{la1}) we have in the Grothendieck
group of the category $\widetilde{\mathcal{O}}_{\bf r'}$
\begin{equation}\label{BtenG}
[B({\bf b}) \otimes e G] =[e \widetilde{M}_{\bf r'}(\epsilon_0)] +
\sum _{i = 1}^{n-1}[e \widetilde{M}_{\bf r'}(\varepsilon_{n-i})[-i
- r_{n-i} + r'_{n-i}]].
\end{equation}

Since $ey^i = y^i e_{n-i}$ it follows that for $1 \leq i \leq n-1$

\[ e \widetilde{M}_{\bf r'}(\varepsilon_{n-i}) = y^i
\mathbb{C}[y^n]\varepsilon_{n-i} \]
 and so
  \begin{equation}\label{eq2}
 p(e \widetilde{M}_{\bf r'}(\varepsilon_{n-i}),s) = s^{i}(1 -
s^{n})^{-1}.
 \end{equation}
In addition we have
 \begin{equation}\label{eq3}
 p(e \widetilde{M}_{\bf r'}(\varepsilon_0), s)  = (1 -
s^{n})^{-1}.
 \end{equation}
 We combine (\ref{eq2}) and (\ref{eq3}) with equation (\ref{BtenG}) and then use equation
 (\ref{kk}) to obtain,
 %since $r_0 = r_n = r'_0 = r'_n = 0,$
 \[ p(B({\bf b}) \otimes eG,s) =  \sum_{i = 0}^{n-1}
  p(e \widetilde{M}_{\bf r'}(\varepsilon_{n-i}),s)s^{-i  -  r_{n-i} + r'_{n-i}} \]
  \[=  \sum_{i = 1}^{n} \frac{ s^{r'_i -r_i}}{(1 - s^{n})} \]
 \[ = \frac{\sum_{i = 1}^{n} s^{n \sum_{j = i}^{n-1}b_j}}{1 -s^n} . \]
 Since this is the same as equation (\ref{obar}), the proof is
 complete.
 \section{\bf Concluding remarks.}
We relate our work to that of Hodges, [H] and Crawley-Boevey and
Holland [CBH]. As in [H] we fix a monic polynomial $v(x) \in
\mathbb{C}[x]$ and let $T(v)$ be the algebra generated by the
elements $h, a, b$ such that
\[ha - ah = a, \quad hb - bh = -b, \quad ba = v(h), \quad ab = v(h-1) . \]
Fix ${\bf k} = (k_1, \ldots, k_{n-1}) \in \mathbb{C}^{n-1}$, and
%et $a_i = (n - i + k_i)/n$ for  $1 \leq i \leq n-1,$
let  ${\bf a} = (a_1, \ldots, a_{n-1}) \in \mathbb{C}^{n-1}$ be as
equation (\ref{400}), and $a_n = 0$.

\begin{lem} \label{hodges}
Set $v(x) = \prod_{i = 1}^n (x - a_i).$ Then $T(v) \cong U_{{\bf
k}}.$
\end{lem}
\noindent {\bf Proof.} The algebra $U_{{\bf k}}$ is generated by
$A = ey^n, B = e(d_{{\bf k}}/n)^n$ and $H = e(\theta + n)/n.$ We
have $HA - AH = A,$ and $HB - BH = -B.$ Using equation (\ref{Ty2})
we see that \[BA  = e\prod_{i = 1}^{n} ((\theta + i - k_{i})/n) =
v(H)\] and similarly $AB = v(H-1).$ Hence there is a surjective
ring homomorphism from $T(v)$ to $U_{{\bf k}}$ sending $a, b, h$
to $A, B, H$ respectively. The map is injective since $T(v)$ and
$U_{{\bf k}}$
have the same associated graded ring.\\

 Now define an action
of the symmetric group $S_{n - 1}$ on $\mathbb{C}^{n-1}$ by
$s(x)_i = x_{s^{-1}i}.$ Define the dot action of $S_{n - 1}$ on
$\mathbb{C}^{n-1}$ by $s\cdot x = s(x + \rho) - \rho.$
\begin{cor} \label{hodges1}
For all $s \in S_{n - 1}, U_{{\bf k}} \cong U_{s\cdot{\bf k} }.$
\end{cor}
%\noindent {\bf Proof.} This is clear since the polynomial $v$ is
%invariant under a permutation of its roots.\\
\indent Now suppose that $v(x) = w(x)u(x)$ is a factorization of
$v$ with $u, w$ monic.  Hodges considers the $T(v)$-module
  \begin{equation}\label{eq300}
 P_w = T(v)a +
T(v)w(h).
 \end{equation}
 The next result is [H, Lemma 2.4].
\begin{thm} \label{hodges2}
The left $T(v)$ module $P_w$ is projective if  $u(x)$ and $w(x)$
are relatively prime, and is a generator if $u(x)$ and $w(x+1)$
are relatively prime.
\end{thm}
We compare $P_w$ to the  $U_{{\bf k}'}-U_{{\bf k}}$-bimodules
$B_p({\bf k}) .$ To do this we now assume that let $v(x) =
\prod_{i = 1}^n (x - a'_i),$ where $a'_i = (n - i + k'_i)/n$ for
$1 \leq i \leq n-1,$ and $a'_n = 0$ so that $T(v) \cong U_{{\bf
k}'}.$  We identify $T(v)$ with $U_{{\bf k}'}$ using this
isomorphism. By Corollary \ref{hodges1} we may assume that $w(x) =
\prod_{i = 1}^p (x - a'_i).$ Then by equation (\ref{Ty2})
$e(d_{{\bf k}'}/p)^ny^n = w(H).$  So comparing equations
(\ref{eq1}) and (\ref{eq300}) we see that $P_w$ is identified with
$B_p({\bf k}).$  Now it is easy to see that Theorem \ref{BC} is
equivalent to Theorem
\ref{hodges2}.\\

For a finite subgroup $\Gamma$ of $SL_2\mathbb{C},$ and $\lambda $
a central element in $\mathbb{C} \Gamma,$ Crawley-Boevey and
Holland define the algebras
\[ S^\lambda = (\mathbb{C}<x, y>\Gamma)/(xy - yx - \lambda),\]
and $\mathcal{O}^\lambda = eS^\lambda e,$ [CBH].  For $\Gamma$
cyclic of order $n,$ we compare the $S^\lambda$ to the algebras
$H_{\bf k}.$ By equation (\ref{Ty1})
  \begin{equation}\label{eq200}
  d_{\bf k} y  - yd_{\bf k}  = 1 +
 \sum_{j = 0}^{n-1} (k_j - k_{j+1})e_j =
 \sum_{j = 0}^{n-1} \lambda_{j}e_j,
 \end{equation}
 where  $\lambda_j = (1/n) + k_j - k_{1+j}.$
Note that the trace of $ \lambda = \sum_{j = 0}^{n-1}
\lambda_{j}e_j$ on the regular representation of $\Gamma$ equals
$1.$ Also if we are given $\lambda \in \mathbb{C} \Gamma$ with
trace $1,$ there is a unique solution to equation (\ref{eq200})
with $k_0 = k_{n} = 0$ and ${\bf k} \in \mathbb{C}^{n-1}.$ Clearly
we have $S^\lambda \cong H_{\bf k},$ and $\mathcal{O}^\lambda
\cong U_{\bf k}.$\\

  \end{document}